\documentstyle[pb-diagram,lamsarrow,pb-lams,amssymb,amsfonts]{amsart}

\def\mystretch{1.2}
\def\baselinestretch{\mystretch}

\def\baselinestretch{\mystretch}

\newtheorem{theorem}{Theorem}[section]

\newtheorem{proposition}[theorem]{Proposition}
\newtheorem{lemma}[theorem]{Lemma}
\newtheorem{corollary}[theorem]{Corollary}

\theoremstyle{definition}
\newtheorem{definition}[theorem]{Definition}

\newtheorem{remark}[theorem]{Remark}

\newtheorem{ack}{Acknowledgments}

\theoremstyle{remark}

\def\abstract{\noindent \vspace{0.3cm}\large{\bf Abstract}\\ 
\small\def\baselinestretch{1}\normalsize}

\newcommand{\W}{{\mathcal {W}}} 
\newcommand{\La}{{\mathcal {L}}} 
\newcommand{\A}{{\mathbb {A}}^\hbar}
\newcommand{\C}{{\mathbb {C}}}
\newcommand{\Z}{\mathbb {Z}}
\newcommand{\R}{\mathbb {R}}
\newcommand{\goth}{\mathfrak}
\newcommand{\U}{{A}}
\newcommand{\F}{{A}}
\newcommand{\disph}{\displaystyle\frac{i}{\hbar}}
\newcommand{\halfh}{\displaystyle\frac{i\hbar}{2}}
\newcommand{\B}{{^{\F}\!\W_B}}
\begin{document}

\title{Deformation Quantization of Symplectic Fibrations.
}
\author {Olga Kravchenko}
\address{IRMA, Universit\'e Louis Pasteur, 7 rue Ren\'e Descartes,
67084,
Strasbourg, France, ok@@alum.mit.edu}	

\maketitle

\vspace{1cm}

\begin{abstract}

\noindent A symplectic fibration is a fibre bundle in the symplectic category. 
We find the relation between deformation quantization of the base and 
the fibre, and the total space.
 
\noindent We use the weak coupling form of Guillemin, Lerman, Sternberg and 
find the characteristic class of deformation of symplectic fibration.
 
\noindent We also prove that the classical moment map 
could be quantized if there
exists an equivariant connection.
 
\noindent Along the way we touch upon the general
question of quantization with values in a bundle of algebras.
 
\noindent We consider Fedosov's construction of deformation 
quantization in general. In Appendix we show how to 
calculate step by step the Fedosov connection, 
flat sections of the Weyl algebra bundle 
corresponding to functions and their star--product. 

\end{abstract}

\vspace{1cm}

\vspace{0.5cm}


\noindent {\bf Key--words:}  deformation quantization, symplectic fibrations, 
moment map, Weyl algebra, Fedosov quantization.

\vspace{0.5cm}

\noindent {\it Address:}
IRMA, Universit\'e Louis Pasteur, 7 rue Ren\'e Descartes,
67084, Strasbourg, France 

\noindent {\it Email:} ok@@alum.mit.edu 	

\newpage

\setcounter{section}{0}
%%%%%%%%%%%%%%%%%%%%%%%%%%%%%%
\section{Introduction: statement of the problem and the main theorem.} 

Quantization is a map from functions on a (phase) space to operators on
some Hilbert space. It involves a parameter (usually the Planck's
constant $h$ or $\hbar = \frac{h}{2\pi i}$).  The product of two
operators is given by some series in $\hbar$.  The inverse of a
quantization map allows one to get a noncommutative product on
functions, namely, by taking the inverse of the product of corresponding
operators. There is also the so called correspondence principle. The
zero degree term in the $\hbar$-decomposition of a noncommutative
function (series in $\hbar$ with functional coefficients) should give a
commutative product of functions. The term of degree one linearly depend
on the Poisson bracket. This allows one to think of quantization as
deformation of a structure of the algebra (Poisson algebra of functions
on a manifold). In the formal deformation quantization one does not
consider questions of convergence of series in $\hbar$. However, the
formal deformation turns out to be a useful tool for describing global
properties of a manifold.  This concept of deformation quantization was
described in \cite{Flato}.  When we say quantization we mean deformation
quantization.

 Fedosov found a beautiful construction of deformation quantization
\cite{F}, \cite{F2} which we use here. Our  article can be considered as
an exercise in his domain, it is in fact a generalization of the first
part of his article \cite{F1} and owes a lot to the methods of this article.

Quantization of a symplectic (or Poisson) manifold $M$ is 
a construction of a noncommutative associative product
on $M$.  It is called a $\ast$--product.  $\ast$--product is defined as 
a product on the
algebra of series in a formal variable $\hbar$ with 
coefficients in $C^\infty (M)$. This noncommutative algebra ${\A}$
should be a deformation
of the algebra of functions on the manifold, $C^\infty (M)$.  
Let $(M, \omega_0)$ be a symplectic manifold. 
Then the symplectic form $\omega_0$ defines a Lie--algebra structure on 
$C^\infty(M)$, called the  Poisson structure.
 (The Jacobi identity follows from closedness of the form $\omega_0$.)
For   $f, g \in C^{\infty}(M)$ let  $\{f, g\} = (df)^\sharp (g) $,
where $\sharp: {\mathcal T}^*M \to {\mathcal T}M$, 
defined by $\omega_0$ (see below
section (\ref{sec:sharp})).
\begin{definition}
Deformation quantization of a symplectic manifold $(M, \omega_0)$ is 
an associative algebra $({\A}, \ast)$ over $\C[[\hbar]]$ with an isomorphism
$\epsilon:  {\A}/ \hbar {\A}  \to C^\infty(M)$ s.t.
\begin{enumerate}
\item ${\A}$ is flat as a $\C[[\hbar]]$--module,
\item ${\A}$ is separable and complete in $\hbar$--topology,
\item for any $f, g \in C^{\infty}(M)$ 
\[
\epsilon(\disph (\check f \ast \check g - 
\check g \ast \check f) )
=\{f, g\}
\]
for $\check f$ and $\check g$ s.t. $\epsilon\{\check {f}\} = f $ and 
$\epsilon\{\check{g}\} = g,$
\item 
The structure of the product $\ast$ on ${\A}/ \hbar^n {\A} $ for all 
$n \geq 0$ is given by bidifferential operators.
\end{enumerate}
\end{definition}

DeWilde--Lecomte \cite{DWL} and also 
Fedosov \cite{F}  proved that on any symplectic manifold there exists a
quantization. 

The following idea lies behind the Fedosov construction (see 
\cite{D}): a
Koszul--type resolution is considered for $C^{\infty}(M)[[\hbar]]$.
Each term of the resolution has a noncommutative algebraic structure
hence providing the  algebra of functions with a new noncommutative product.
Fedosov constructs  such a resolution by using the differential 
forms on the manifold with values in the  Weyl--algebra bundle. 
The main step then is to find  a differential on it which respects the 
algebra structure. This
differential is called Fedosov connection and is obtained by an iteration
procedure from a torsion free 
symplectic connection on the  manifold. 

Lichnerowitz \cite{L} showed that any connection on a symplectic
manifold gives rise to a torsion--free symplectic connection. Hence one
can get the Fedosov connection from any connection on a tangent bundle:
first adding some tensor to make it symplectic (section
\ref{sec:sharp}) and then applying the iteration procedure.

We introduce a notion of an $F$--manifold:
\begin{definition}
An $F$--manifold  $\Phi$ is the following triple 
(manifold, deformation of a symplectic form,
a connection) : 
\[
 \Phi = (M, \omega, \nabla),
\]
 where $\omega = \omega_0 + \hbar \alpha $, and 
$\alpha \in \Gamma (M, \Lambda^2 T^*M) [[\hbar]]$, a series
in $\hbar$ with coefficients being closed $2$--forms on the manifold.
\end{definition}
 In a recent article \cite{GRS} a similar object is
called a Fedosov manifold, namely, a symplectic manifold together with a
symplectic connection. Indeed,
these three objects $(M, \omega_0, \nabla)$ define 
the first three terms in the $\ast$--product: 
\begin{itemize}
 \item Classically, a manifold can be considered as an algebra of
functions on it, that is $M$ defines the structure of the commutative
product.
 \item A symplectic form defines the Poisson structure and hence the 
term at $\hbar$.
 \item A connection defines the term at $\hbar^2$ (as follows from \cite{L}).
\end{itemize}
 It turns out  that these three terms determine the higher terms. The 
deformation quantization  theorem (\cite{F})
can be stated as follows
\begin{theorem}(Fedosov)
\label{th:fedosov}
An $F$--manifold  $(M, \omega, \nabla)$  uniquely determines 
a $\ast$-product on the underlying manifold $M$.
 \end{theorem}
Deligne \cite{De} and  Nest and Tsygan \cite{NT}  showed
 that the class of 
isomorphisms of quantizations of a symplectic manifold $M$ 
is  determined by the class of the form 
$\omega$ in $ H^2(M)[[\hbar]].$ It is called a 
 {\em characteristic class of deformation}.

\vspace{5mm}
 We study the deformation of the {\em twisted products}\/ of
two $F$--manifolds $(B, \omega_B, \nabla^B)$ and $(F, \sigma, \nabla^F)$. 
Our  question is the following: how to define the product of two
$F$--manifolds and what the $\ast$--product on the total space is,
that is 
what a ``twisted product of quantizations'' is.

 We show under certain assumptions that a twisted product of
two $F$--manifolds is again an $F$--manifold. So we want to relate
the $\ast$--products on these three manifolds.

One can regard the product $ M = B \ltimes F$ as a
fibre bundle $M \to B$ over a symplectic base $B$
with a symplectic  fibre $F.$ 
Obviously, the product depends on how twisted the symplectic fibration 
$ M = B \ltimes F$ is.
This can be described by a  connection on $ M \to B$ which should be   
compatible with the symplectic form on $M$.

 When $G$, the structure group  of the bundle $M \to B$, acts by
symplectomorphisms on fibres this bundle  is a symplectic fibration 
(see \cite{GLS}, see bellow section \ref{sec:sf}). 
The total space  $M$ is symplectic and the  fibres 
are symplectic submanifolds of it.
 We need  the action to be  Hamiltonian and 
$G$ to be a finite dimensional Lie group.

 The case of $B$ being a cotangent bundle of some differentiable
manifold $X$ was first developed by Sternberg \cite{S}.  His
construction describes the movement of a "classical particle" in
Yang--Mills field for any gauge group $G$ and any differentiable $X$.
 Weinstein  gave a general construction of a symplectic form
on $M$ in \cite{W}.

Here we are mainly interested in the particular case of
the quantization with coefficients in the auxiliary bundle of 
fibrewise quantizations. 
This leads us to a more general case of quantization with coefficients
in some bundle of algebras (some examples of such bundles are considered
in \cite{F2}). Our case is  somewhat different
 from the quantization with
values in most  other auxiliary bundle and is more difficult to
perform. Fibres of the auxiliary bundle
obtained from the symplectic fibration structure 
are  noncommutative algebras. It is this noncommutativity of fibres
which  makes the quantization procedure more complicated.

However we believe that it is useful to understand its
mechanism in order to see that quantization is a fundamental notion like
some homology theory and hence it should respect a fibre bundle structure.

We make a new definition:
\begin{definition}
An $F$--bundle (with an underlying manifold $B$) is a triple 
\[
\Psi = (\Phi,  \U, \nabla ^{\F})
\]
where  
\begin{itemize}
\item
$\Phi =(B, \omega^B, \nabla^B)$ is an $F$--manifold, 
\item $\U$ is  an 
auxiliary bundle of algebras over $B$,
\item $\nabla ^{\F}$ is  a covariant 
derivative on $\U$, which respects the algebraic 
structure on the fibres.
\end{itemize}
\end{definition}

We construct an $F$--bundle from a symplectic fibration $M \to B$.
Each fibre of $M \to B$ is an $F$--manifold, so we can quantize the fibres. 
If the connection 
$\nabla^F$ is
$G$--invariant we can consider a new bundle $\U$
over the base $B$. This bundle $\U$ is a bundle of algebras
of quantized functions on fibres. The fibre of $\U$ over a point
$b$ is the algebra of quantized functions on the fibre of $M \to B$ at
the point $b$:
 \begin{equation}
 {\U}_b = {\A} (M_b).
        				\label{eq:aux} 
\end{equation}

The bundle $\U$ is defined by (\ref{eq:aux}).
The
covariant derivative $\nabla^{\F}$ on the bundle $\U$
is determined by the connection on the bundle $ M \to B$. The
construction of the covariant derivative is carried out in section
(\ref{sec:derivative}) so that it respects
the algebraic structure, that is 
 it satisfies Leibnitz rule with respect to the product on $\U$.

Such an $F$--bundle corresponds to an $F$--manifold modeled on $M$, the
total space of the symplectic fibration, and hence provides a
quantization of the total space. Our main theorem is:
 \begin{theorem}
\label{th:main}
 Consider  a symplectic fibration $M \to B$ with 
a standard fibre being an $F$--manifold
$(F, \sigma, \nabla^F$ and the base an $F$--manifold 
$\Phi = (B, \omega^B, \nabla^B)$. 
 An  $F$--bundle $ (\Phi, \U, \nabla ^{\F}) $, 
where  ${\U}_b = {\A} (M_b)$ 
gives a quantization of the underlying 
manifold $B$ with values in the auxiliary bundle. 
This also defines the quantization of the total space $M$.
\end{theorem}

The main claim of this theorem is that quantization of the base with
values in the auxiliary bundle (\ref{eq:aux}) corresponds to a certain
$F$--manifold $\Phi = (M, \omega, \nabla)$, where $\omega$ is
 a polynomial in $\hbar$ starting from a symplectic form on $M$.

To carry out the program first of all one has to construct a symplectic
form on $M$.  There is a one--parameter family of symplectic forms on
the total space. The construction involves the notion of weak coupling
limit of Guillemin, Lerman and Sternberg \cite{GLS}.  The behavior of
the $\ast$--product when this parameter tends to zero gives us a way to
understand the relation between quantizations of the base and the fibre
with the quantization of the total space.

Our main Theorem (\ref{th:main}) can be reformulated as a statement 
about the solutions of two 
equations given in Theorem (\ref{th:main1}). 
We construct a twisted Weyl algebra bundle corresponding to 
$F$--bundle. Then we prove that 
\begin{enumerate}
\item There exists $r$, a $1$-form on the base with values in the twisted 
Weyl bundle, such that the initial connection becomes flat 
when one adds $r$ to it.
\item For each section ${\mathbf a}$
of the auxiliary bundle there exists only one corresponding 
flat section of the twisted Weyl bundle.
\end{enumerate}

Fedosov's quantization procedure is discussed in Section
\ref{chap:fedosov}, calculations and examples are given in the Appendix.
The classical setup for symplectic fibrations is discussed in Section
\ref{chap:sf} and the quantization of the moment map is presented in Section
\ref{qmm}.  Main results about the quantization of symplectic fibrations
are given in Section \ref{chap:main}, examples are discussed in the last
section.

\textbf{Notations.} Repeated indices assume summation.  

Grading and filtration of the Weyl algebra bundle 
are $\Z$--grading and  $\Z$--filtration, we  do not use the natural
${\Z}_2$--grading on the differential forms.

Let ${\mathcal E}$ be a bundle over some manifold $M$. Then
${\mathcal A}^n (M, {\mathcal E})$ denotes  $C^\infty$--sections 
of $n$--form bundle with values in the bundle $\mathcal E$,
\[
{\cal A}^k (M,{\cal E})= {\Gamma} 
(B, \Lambda^k {\mathcal T}^*B \otimes {\cal E}). 
\]
${\mathcal A}^n (M)$ denote the bundle of $n$--forms on $M$, and 
\[
{\mathcal A} (M, {\mathcal E})= 
\oplus_{n=0}^\infty{\mathcal A}^n (M, {\mathcal E}).
\]

The term ``connection'' is used in two senses, for a covariant
derivative on any vector bundle, usually denoted by $\nabla$ and also
for a connection on a fibre bundle i.e. a splitting of the tangent
bundle to the total space of a fibre bundle into a sum of  a vertical  
and a horizontal subbundles.

\begin{ack} I am deeply grateful 
to Ezra Getzler, Richard Melrose, Boris Tsygan and 
David Vogan for help and encouragement. My thanks also go to 
Alex Astashkevich,
Paul Bressler, David Ellwood, Pavel
Etingof, Daniel Grieser, Dima Kaledin, Eugene Lerman, and  Andras Szenes 
for many fruitful discussions.

 This work started at
mathematics department at MIT and completed during my stay at IHES
as a European postdoctoral fellow. It is a pleasure to acknowledge
the support of these two institutions.
 \end{ack} 

%%%%%%%%%%%%%%%%%%%%%%%%%%%%%%%

\section{Generalities on deformation quantization.}
\label{chap:fedosov}

The subject of this section becomes nowadays fairly standard (see for
example an excellent introduction to Fedosov quantization \cite{Le}).

\subsection{Weyl algebra of a vector
space.\\}
 
Let $E$ be a vector space with a non--degenerate skew-symmetric form $\omega$.
  The algebra of polynomials on $E$ is the algebra of symmetric
powers of $E^*$, $S(E^*)$, and it has a skew-symmetric form on it 
which is dual to $\omega$.
Let $e$ be a point in $E$ and $\{e^i\}$ denote its linear coordinates in
$E$ with respect to some fixed basis.  Then $\{e^i\}$ define a basis in $E^*$. 
Let $\omega^{ij}$ be the
matrix for the skew-symmetric form on $E^*$. Let us consider the power
series in $\hbar$ with values in $S(E^*)$: 

\begin{definition}
The Weyl algebra $W(E^*)$ of a vector space $E^*$ 
is an associative algebra
 \[ 
W(E^*) = S(E^*) [[\hbar]]: \ \ \ 
 a(e,\hbar) = \sum_{k \geq 0} a_k(e) \hbar^k,
\] 
 with  the product structure given by the Moyal--Vey product:
	\begin{equation}
	\left.
 a\circ b (e,\hbar)= \exp \  \{- \halfh \omega ^{kl} \frac{\partial  }
{\partial
 x^k} \frac{\partial  }{\partial  z^l}\} \ \  a(x, \hbar) \ 
b(z, \hbar)\  \right| _{x=z=e}
							 \label{eq:Moyal}
	\end{equation}
 \end{definition}

\noindent 
The Lie bracket is defined with respect to this product.
We can look at this algebra as at a completion of 
the universal enveloping algebra of the 
Heisenberg algebra on $E^* \oplus \hbar \C $, namely, the algebra 
with relations
\[
e^i \circ e^j - e^j\circ e^j = - i {\hbar}\omega^{ij}
\]
where $\omega ^{ij} = \omega (e^i, e^j)$ defines  a  Poisson bracket
on $E^*$.

Hence one can define the Weyl algebra as
\[
W = U (E^* \oplus \hbar {\C}).
\]

Let us consider the product of the Weyl algebra and the 
exterior algebra of the space $E^*$: $W(E^*)\otimes \Lambda E^*$.
Let $dx^i$ be the basis in $\Lambda E^*$ corresponding to $e^i$ in $W(E^*)$.
 
There is a decreasing  filtration on the Weyl 
algebra $W(E^*)$: $
W_0 \supset W_1 \supset W_2 \supset ...$  given by the degree of
generators. 
$e$'s have degree 1 and $\hbar$ has degree 2: 
\[
W_p = \{ \text{elements with degree}  \geq p \}. 
\]

One can define a grading on
$W$ as  follows
\[
gr_i W = \{ \text{elements with degree}  = i \}. 
\]
it is isomorphic to 
$W_i/W_{i+1}$.
One can see  that  the product (\ref{eq:Moyal}) preserves the grading.

\begin{definition}
\label{def:gr}
An operator on  $W(E^*)\otimes \Lambda E^*$ 
is said to be of degree $k$ if it maps 
$W_i\otimes\Lambda E^*$ to $W_{i+k}\otimes\Lambda E^*$ for all $i$.
\end{definition}

\noindent 
Such an operator defines maps $ gr_i W\otimes\Lambda E^*$ to
$gr_{i+k}W\otimes\Lambda E^*$ for all $i$.

\begin{definition} 
Derivation on  $W(E^*)\otimes \Lambda E^*$ is 
a linear operator which satisfies the Leibnitz rule:
\[
D (ab) = (Da) b + (-1)^{\tilde{ a}\tilde {D}}a (Db)
\]
\end{definition}
\noindent 
where $\tilde{a }$ and $\tilde{D }$ are  corresponding degrees. 
It turns out that all linear derivations are inner operators.
\begin{lemma}
 \label{l:inner}
 Any linear derivation $D$ 
on $W (E^*) \otimes \Lambda E^*$ is inner, i.e.
 namely there exists such $v\in W(E^*)$ so that 
$D a = \disph [v,a]$ for any $a \in W (E^*)$
\end{lemma}
\begin{pf} 
Indeed,
$\frac{\partial}{\partial e^i} a = \frac{i}{2 \hbar}[ \omega_{ij} e^j, a].$\\  
So for any derivation one can get a formula:
$D a = \disph [ \frac{1}{2} \omega_{ij} e^i D e^j , a].$  
\end{pf}
One can define two natural  operators on the algebra 
$W(E^*)\otimes \Lambda E^*$: $\delta$ and 
$\delta^*$ of degree $-1$ and $1$ correspondingly.
$\delta$ is the lift of the ``identity'' operator 
\[
u: e^i \otimes 1 \to 1 \otimes dx^i
\]
and  $\delta^*$ is the lift of its inverse.
On monomials  $e^{i_1}\otimes...\otimes e^{i_m}\otimes dx^{j_1} \wedge ...
\wedge dx^{j_n} \in W^m(E^*)\otimes \Lambda^n E^*$ 
$\delta$ and 
$\delta^*$ can be written as follows:
\begin{eqnarray*}
\lefteqn{
\delta: \ e^{i_1}\otimes...\otimes e^{i_m}
\otimes dx^{j_1} \wedge ...\wedge dx^{j_n}
\mapsto}\\
& &\sum _{k=1}^m
e^{i_1}\otimes...\widehat{e^{i_k}}...\otimes
e^{i_m}\otimes dx^{i_k} \wedge dx^{j_1} \wedge ...\wedge dx^{j_n}\\
\lefteqn
{\delta^*: \  e^{i_1}\otimes...\otimes e^{i_m}
\otimes dx^{j_1} \wedge ...dx^{j_n}
\mapsto}\\
& &\sum _{l=1}^n (-1)^l e^{j_l} \otimes
e^{i_1}\otimes...\otimes e^{i_m}\otimes dx^{j_1} \wedge ...
\widehat{dx^{j_l}} ... \wedge dx^{j_n}.
\end{eqnarray*}
Let  $a_0$ be  a projection  of $a \in W(E^*)\otimes \Lambda E^*$ to 
$ gr_0W(E^*)\otimes \Lambda^0 E^*$, which is the center of the algebra, i.e.
the summands in  $a$ which  do not contain neither  $e$-s or $dx$-s.
\begin{lemma}
\label{l:delta}
Operators $\delta$ and $\delta^*$ have the following properties:
\[
\delta a = dx^j  \frac{\partial a }{\partial  e^j} =
\ [ -\frac{i}{\hbar} \omega_{kl} \ e^k dx^l, \ a], \quad
\delta ^* a = y^j \iota_{\frac{\partial }{\partial  x^j}} a, \quad
\delta^2 = {\delta^*}^2 = 0
\]
On monomials from $gr_m W(E^*)\otimes \Lambda^n E^*$ 
\[
\delta \delta^* + \delta^* \delta = (m+n) Id,
\] 
where $Id$ is the identity 
operator.
Any element $a \in gr_m W(E^*)\otimes \Lambda^n E^*$ has a decomposition:
\[
a = \frac{1}{m+n}(\delta \delta^* a + \delta^* \delta a)  + a_0.
\]
\end{lemma}

%%%%%%%%%%%%%%%%%%%%%%%%%%%%%%%%

\subsection{Symplectic connections (covariant derivatives).\\}

\label{sec:sharp}

The term  {\em symplectic connection} in this section in fact must be
changed to {\em symplectic covariant derivative} to avoid confusion with
another symplectic connection notion in the next Chapter. However there
is already an  established practice to call a covariant derivative a 
connection which we decided to follow here. We hope that one can get used to 
distinguish one from the other from the context.

Let us consider connections on a manifold $M$.
	\begin{proposition}
Let $\omega$ be a skew--symmetric 2--form on ${\mathcal T}M$. 
Then $\omega$ must be 
closed in order for
torsion--free connection $\nabla$ preserving this form to exist.
	\end{proposition}
\begin{pf} The skew--symmetry of $\omega$ is the following condition: 
$\omega (X,Y) = - \omega (Y,X)$.  
The connection $\nabla$ is torsion--free when 
$\nabla_X Y -\nabla_Y X  = [X,Y]$. 
Suppose such $\nabla$ exists. Then it preserves the form
$\omega$ when $\nabla \omega = 0$. This means that for all $X,Y, Z \in {\mathcal T}M$:
\begin{equation}
\nabla_X (\omega ( Y, Z)) = \omega (\nabla_X Y , Z) +\omega ( Y ,\nabla _X Z) 
				\label{eq:nabla}
\end{equation}
Since $ \omega ( Y, Z)$ is a function 
 $\nabla_X ( \omega ( Y, Z)) = X \omega  ( Y, Z)$.
Then, 
\begin{eqnarray*}
X \omega  ( Y, Z) & - & Y \omega  ( X, Z) \ \ + \ \  Z \omega ( X, Y) \\
 & = & \omega (\nabla _X Y , Z) - \omega (\nabla _X Z, Y) -
 \omega (\nabla _Y X, Z) \\
 &  &\mbox{}+ \omega (\nabla _ Y Z, X) +
 \omega (\nabla _Z X, Y) - \omega (\nabla _Z Y, X)\\
& = & \omega ([X, Y], Z) - \omega ([X,Z], Y) + \omega ([Y, Z], X)
\end{eqnarray*}
which is exactly the condition $d \omega = 0$.
\end{pf}

Notice, that in the Riemannian case when the form is symmetric, there is
a unique torsion free connection compatible with the form, the
Levi--Civita connection. In the case of a skew--symmetric form there
are plenty of connections compatible with the form, provided that the
form is closed. So the statement of uniqueness of Levi--Civita
connection in the Riemannian case is substituted by the 
requirement for the form to be closed in the skew--symmetric setting.

Here we are mostly interested in the case when $M$ is a 
symplectic manifold, i.e. there is a symplectic form  $\omega$ 
(a closed and nondegenerate $2$-form  on ${\mathcal T}M$).
\begin{definition}
\label{def:con}
A connection which preserves  a symplectic form is called a 
symplectic connection.
\end{definition}
Any connection on a symplectic manifold gives rise 
 to a symplectic connection:

\begin{proposition} \cite{L} \cite{MRR}.
Let $\omega $ be a closed nondegenerate 2--form. Then for every 
connection $\nabla$ there exists a three--tensor $S$, such that 
\[
\tilde{\nabla} = \nabla + S
\]
is a connection on ${\mathcal T}M$ compatible with $\omega$.\\ 
Then for $X,Y \in {\mathcal T}M$
\[
\hat{\nabla}_X Y = \tilde{\nabla}_X Y - \frac{1}{2} Tor (X,Y)
\]
defines a torsion--free connection compatible with the form $\omega$. Here 
$2$--form $Tor$ is the  torsion of $\tilde{\nabla}$
\[
Tor (X,Y) = \tilde{\nabla} _X Y -   \tilde{\nabla} _Y X -
  \tilde{\nabla} _{[X,Y]}
\]
\end{proposition}
Then $S$ is defined as follows:
\[
S_X Y = \frac{1}{2} \{ (\nabla_X \omega )( Y, .) \}^{\sharp},
\]
where $\sharp:  {\mathcal T}^*M \to {\mathcal T}M $ is  the inverse  
to $\flat$, given by :
\begin{eqnarray*}
\flat:{\mathcal T}M \to {\mathcal T}^*M\\
u^\flat = \omega (u, .) \ for \ \ u \in {\mathcal T}M.
\end{eqnarray*}
% Let ${\cal A}^* (M)$ denote the algebra of differential forms on $M$.
Symplectic connections form an affine space with the  associated vector 
space    ${\mathcal A}^1(M,sp(2n)),$
Lie algebra  $sp(2n)$ valued one--forms on $M$.

\subsection{Deformation quantization of a symplectic manifold.\\}
%%%%%%%%%%%%%%%%%%%%%%%%%%%%%%%%%%%
Let  $M^{2n}$ be a symplectic manifold with a symplectic form 
$\omega$. In local 
coordinates at a  point $x$:
      \[
\omega = \omega_{ij} d x^i \wedge dx^j.
\]
The symplectic form on a manifold $M$ defines a Poisson  
bracket on functions on $M$.
For any two functions $u, v \in C^{\infty}(M)$:
	\begin{equation}
	\{u,v \} = \omega^{ij} \frac{\partial  u}{\partial  x^i} 
\frac{\partial  v}{\partial  x^j}    \label{eq:poisson}
	\end{equation}
where $(\omega^{ij})=(\omega_{ij})^{-1}$. 

We can define the bundle of Weyl algebras  ${\cal W}_M$, with 
the fibre at a point $x \in M$
being the Weyl algebra of ${\mathcal T}^*_x M$. Let  
$\{e^1, ... e^{2n}\}$ be $2n$ 
generators in  ${\mathcal T}^*_x M$, corresponding to $dx^i$. The form $\omega^{ij}$ defines 
pointwise Moyal--Vey product. 

\noindent

The filtration and the grading in ${\cal W}_M$ are inherited 
from $W({\mathcal T}^*_xM)$ at each point $x \in M$. 
Denote by  ${\cal W}^i$ the $i$-th graded component in ${\cal W}_M$: 
 \[
{\cal W}_M=\oplus_i {\cal W}^i
\]
 A symplectic connection, $\nabla$, satisfying (\ref{eq:nabla})
can be naturally lifted to act on any symmetric power 
of the cotangent bundle (by the Leibnitz rule) and since 
the cotangent bundle ${\mathcal T}^*M \cong {\cal W}^1$ we 
can lift $\nabla$ to be an operator on sections 
$\Gamma (M,{\cal W}^i)$ with values in 
$\Gamma (M, {\cal W}^i \otimes {\mathcal T}^*M)$.
By abuse of notations this 
operator is also called $\nabla$. 

\noindent It preserves the grading, in other words 
it is an operator of degree 
zero.
It is clear that in general this connection is not flat: $\nabla ^2 \not= 0$.
Fedosov's idea is that for ${\W}_M$ bundle one can add to 
the initial symplectic connection some 
operators not preserving the  grading
so that the sum gives a flat connection on the Weyl bundle.

\begin{theorem}(Fedosov.)
\label{Fedosov}
There is a unique set of operators 
$r_k : \Gamma (M,{\cal W}^i) \to
\Gamma (M, {\mathcal T}^*M \otimes {\cal W}^{i+k})$ such that 
\begin{equation}
	D = - \delta + \nabla + r_1 + r_2 + \ldots 
				  \label{eq:connection}
	\end{equation}
is a flat connection and 
\[
\delta^* r_i =0.
\]
There is a one-to-one correspondence between formal series in $\hbar$ 
with coefficients in smooth functions $C^\infty(M)$
and  horizontal sections of this connection:
\begin{equation}
\label{eq:Q}
Q: C^\infty(M)[[\hbar]] \to  \Gamma_{\text{flat}} (M,{\cal W}_M).
\end{equation}
\end{theorem}
\noindent Main idea of the proof is to use the following complex:
 \begin{equation}
\label{comp:del}
0 \stackrel{ } \to \Gamma (M,{\W} ) 
\stackrel{\delta} \to {\mathcal A}^1(M,{\W} ) 
\stackrel{\delta} \to {\mathcal A}^2(M,{\W} ) \stackrel{\delta} \to \ldots. 
 \end{equation}
This complex is exact since $\delta$ is homotopic to identity  by $\delta^*$. 
An equation for $r_i$ for each $i > 1$ has the form 
 \begin{equation}
 \label{eq:r}
\delta (r_i) = \text{function}(\nabla, r_1, \ldots, r_{i-1}).
 \end{equation}
However it is not difficult to show that  this function is in the kernel
of $\delta$ hence $r_i$ exists.

First few steps in the construction of $D$, its flat sections and
$\ast$--product in coordinates are given in the Appendix.

The  noncommutative structure on the Weyl bundle
determines a $\ast$--product
on functions by  this correspondence, namely for two functions 
$f,g \in
\C^\infty(M)[[\hbar]]$;
 \begin{equation}
\label{eq:ast}
f\ast g = Q^{-1} (Q(f) \circ Q(g)). 
\end{equation}

In fact, the equation $D^2 = 0$
is just the Maurer--Cartan equation for a flat connection. 
One can see the analogy with 
Kazhdan connection \cite{GKF} on the algebra of formal vector fields.
Notice that  $\delta = dx^i
 \frac{\partial }{\partial e^i}$  is of degree $-1$. The flatness of
the connection is given by the  recurrent  procedure, namely starting from the 
terms of degree $-1$ and $0$ one can get other terms step by step. 
 While Kazhdan connection does not have a parameter involved it has the
same structure -- it starts with known $-1$ and $0$ degree terms. Other
terms are of higher degree and can be recovered one by one.

Let us also
mention  here that the connection $D$ can be written as a sum of two terms
-- one is a derivation along the manifold, the usual differential $d$, and
the other is an endomorphism of a fibre, let us call it $\Gamma $. 
Since all endomorphisms
are inner one can write it as an adjoint action with respect to the Moyal 
product. $\Gamma $ acts 
adjointly by an 
operator from $\Gamma (M,{\cal W})$ to $\Gamma (M, {\mathcal T}^*M
\otimes {\cal W})$. 
 \begin{equation}
\label{eq:con}
D = d + \Gamma = d + \disph[{\gamma}, \cdot]_\circ,
 \end{equation}
where  ${\gamma} \in \Gamma (M, {\mathcal T}^*M \otimes {\cal W})$.
 Then the  equation $D^2 = 0$ becomes:
 \[
d \Gamma + \frac{1}{2} [\Gamma,\Gamma ]_\circ =0.
\]
The same equation for  ${\gamma}$ then is as follows:
 \begin{equation} 
\label{eq:eq}
\omega + d \gamma + \disph \frac{[\gamma,\gamma]_\circ}{2} = 0,
 \end{equation}
where $\omega$ is a central 2--form.
 This equation states that $D^2$ is given by an
adjoint action of a central  element, so it is zero.
 However it turns out to be very important which exactly form  $\omega$ is 
given in the center by the connection $D$.
 Inner automorphisms of the Weyl algebra are given by the 
adjoint action by  elements of the algebra (Lemma \ref{l:inner}).
Its central extension gives the whole algebra. 
 Curvature of Fedosov connection is zero, however its lift to the
central extension is nonzero and gives the characteristic class
of quantization.
\begin{definition}
\label{def:char}
The characteristic class of the deformation is the cohomology 
class of the form  $[\omega] \in H^2  (M)[[\hbar]].$
\end{definition}

%%%%%%%%%%%%%%%%%%%%%%%%%%%%%%%%%%%%%%%%%%%%%
\noindent {\bf Quantization with values in a bundle of algebras.} 
 This subject
was discussed at length in the book \cite{F1}, but here we want to look
at it from a slightly different angle.
  Given an $F$--manifold $(B, \omega, \nabla^B)$ we 
know how to construct a map $C^\infty (B) \to \Gamma(B, \W_B)$.
Let $\La \to B$ be a bundle of $\C[[\hbar]]$--algebras.
Now we want to generalize the problem of quantization and obtain a map:
 \[ 
Q^{\La}:  \Gamma(B, \La) \to \Gamma(B, \W_B \otimes_{\C[[\hbar]]} \La).
 \]
 In order to do that we need a connection $\nabla^{\La}$ on $\La$. Let
$R^{\La} \in {\mathcal A}^2 (B, \La)$ be its curvature. Also the sections of $\La$
must commute with sections of the Weyl algebra bundle $\W_B$. Then we
can define a connection on $\W_B \otimes_{\C[[\hbar]]} \La$ as the sum
of connections.
\[
  \nabla = \nabla^B \otimes 1 + 1 \otimes \nabla^\La
\]
with the curvature 
\[
 R = R^B \otimes 1 + 1 \otimes  R^{\La}.
\]

%\hspace{-2cm}{\tt Module without a torsion.}
 We define a grading on 
 $\W_B \otimes_{\C[[\hbar]]} \La$ as on  $\W_B$. 
The operator $R^{\La}$ could have a degree if it changes the power of 
 $\hbar$, and since the degree of  $\hbar$ is $2$ it could  only be even:
 $R^{\La} = \sum R^{\La}_{2k}, \ k \leq 0$ where 
\[
R^{\La}_{2k}: {\mathcal A}^q (B, gr_{\cdot}(\La) ) 
\rightarrow  
{\mathcal A}^{q+2}(B,gr_{\cdot +2k}(\La)).
\]
One should add
an extra term in the equation on 
$r_i$ (\ref{eq:r}) for even $i = 2k$ one should add
an extra term
 \[
 \delta (r_{2k}) = \text{function}(\nabla, r_1, \ldots, r_{{2k}-1})
+ R^{\La}_{2k}. 
 \] 
 Then like in Theorem \ref{Fedosov} we can consider a  complex similar to 
(\ref{comp:del})
\[
0 \stackrel{} \rightarrow \Gamma (B,\W_B \otimes_{\C[[\hbar]]} \La ) 
\stackrel{\delta} \rightarrow  
{\mathcal A}^1(B,{\W_B} \otimes_{\C[[\hbar]]} \La ) \stackrel{\delta} 
\rightarrow {\mathcal A}^2(B,\W_B \otimes_{\C[[\hbar]]} \La )
  \stackrel{\delta} \rightarrow \ldots 
\]
with $\delta$ acting only in $\W_B$. This complex is still exact (The  
$ \otimes$--product is over a field ${\C[[\hbar]]}$).
Since  $\delta R^{\La} = 0$ one can find a preimage of $ R^{\La}$
and the reasoning is exactly as before.
 The flat connection and the corresponding flat
sections are constructed  similarly  to the way 
it  is described in Appendix. 

 However in the case when ${\La}$ has a Lie algebra structure and $
R^{\La}$ is an inner action of the form $ R^{\La} = \disph \text{ad} H$
for some $H \in {\cal A}^2(B, \La)$ the procedure changes! The adjoint
action might start from the degree $ - 2$ term and one has to change not
only the equations (\ref{eq:r}), but also the initial $\delta$ to
balance it.  This is exactly what happens in our case of symplectic
fibrations and what makes it more interesting.

%%%%%%%%%%%%%%%%%%%%%%%%%

\section{Symplectic
forms on symplectic  fibrations.\\} 

\label{chap:sf}

\indent In this section we collect the facts known about symplectic
fibrations: we give a definition and construct a one-parameter family of
symplectic forms on the total space.  One can find a nice exposition in
the sixth chapter of the book \cite{McDS}, see also \cite{GLS}. 

The meaning of the term {\em symplectic connection} used in this section
is different from the definition \ref{def:con}. It is a connection which
preserves the symplectic structure on fibres while in the previous
Chapter \ref{chap:fedosov}  the symplectic connection was in fact
a symplectic covariant derivative preserving a symplectic form on a 
symplectic manifold.

\noindent{Symplectic fibrations.}
 \label{sec:sf}
\begin{definition} A symplectic fibration is a locally trivial 
fibration $ \pi: M \to B$ with a symplectic fibre $(F, \sigma)$  
whose structure group preserves the 
symplectic form $\sigma$ on $F$. This means that 
\begin{enumerate}
\item
There is an open cover ${U_\alpha}$ of $B$ and a 
collection of diffeomorphisms 
$
\phi_\alpha: \pi ^{-1} U_\alpha \to
 U_\alpha \times F
$
such that the following diagram commutes: 
\[
   \dgARROWLENGTH=0.6\dgARROWLENGTH
   \begin{diagram}
                  \node{\pi ^{-1}U_\alpha}\arrow[2]{e,t}{\phi_\alpha}
		\arrow{se,b}{\pi}			
 \node[2]{U_\alpha \times F}\arrow{sw,b}{pr}	\\
\node{} \arrow{e,!}      \node{U_\alpha}
   \end{diagram}
\]
\item
For the fibre over $b\in B$, 
$F_b = \pi^{-1}(b)$, let $\phi_\alpha(b)$
denote the restriction of $\phi_\alpha$ to $F_b$ followed by projection
onto $F$, 
$\phi_\alpha(b): F_b \to F$. Then
\[
\phi_{\beta \alpha} (b) =  \phi_\beta(b)\circ\phi_\alpha(b) ^{-1}
\ \ \ \in Symp(F, \sigma)
\]
for all ${\alpha}$, $\beta$ and $b\in U_\alpha \bigcap U_\beta$.
\end{enumerate}
\end{definition}

If  $ \pi: M \to B$ is a symplectic fibration then 
each fibre $F_b$
carries a symplectic structure $\sigma_b \in \Omega^2 (F_b)$ defined by 
\[
\sigma_b = \phi_\alpha (b)^* \sigma
\]
for $b\in U_\alpha$. The form is independent of $\alpha$ as follows from the
definition. Also, if there is a $G$--invariant 
symplectic torsion--free connection 
$\nabla^F$ on $F$ it defines a symplectic torsion--free connection $\nabla_b$
on each fibre $F_b$.
\begin{definition} A symplectic form $\omega$ 
on the total space $M$ of a symplectic fibration 
is called compatible with the fibration $\pi$ if each 
fibre $(F_b,\sigma_b)$ is a symplectic submanifold of 
$(M, \omega)$, with $\sigma_b$ 
being the restriction of $\omega$ to  $F_b$.
\end{definition}

\noindent {\bf Symplectic connections.} Each symplectic 
form compatible with a symplectic fibration $\pi: M \to B$ 
defines a connection on it, i.e. a choice of splitting 
of the following short exact sequence of vector bundles:
\[
0 \to {\mathcal V}\!M \to {\mathcal T}\!M \to \pi^* {\mathcal T}B \to 0.
\]
 Here ${\mathcal V}M$ is the canonically defined bundle of vertical
tangent vectors, i.e. those fields which vanish on functions coming from the
base. In other words the connection is the splitting of the tangent
bundle into the sum
 \begin{equation}
	\Gamma : \quad {\mathcal T}\!M = {\mathcal H}M 
\oplus {\mathcal V}\!M
				  \label{eq:split}
	\end{equation}
such that $ \pi^* {\mathcal T}B = {\mathcal H}M $.
The connection (\ref{eq:split}) is compatible with the
symplectic form if  at each point $x \in M$: 
\[ 
{\mathcal H}_xM : = \{X \in {\mathcal T}\!_xM: \quad 
 \Omega (X,V) =0 \quad  \text{for  all}
\ {V} \in {\mathcal V}\!_x M\}.
\] 
So each symplectic form whose restriction on fibres is nondegenerate
defines a compatible connection. Namely, 
the  horizontal subbundle consists of 
all  vector fields which are perpendicular to the vertical ones
with respect to the symplectic form.

\noindent {\bf Ingredients: a connection on a 
principal bundle and a Hamiltonian action along the fibres.}
Symplectic fibrations are associated fibre bundles 
to the principal bundles with a structure group being the group of
symplectomorphisms of the fibre, so we have a principal $G$--bundle and 
a symplectic manifold $(F,\sigma)$ to start with.
 
Let us first consider a principal $G$--bundle, i.e. a smooth manifold
$P$ with a smooth action $P\times G \longrightarrow P$ which is free and
transitive. Then the quotient $P/G = B$ is a manifold.

 For the principal bundle a connection can be define by a so called
connection 1-form. Namely, the fibres of a vertical subbundle ${\mathcal
V}\!P$ are naturally identified with $\goth g$ under the map: ${\goth g}
\to Vect(P)$ given by the infinitesimal action of $G$ on $P$.
 \[
 X \in {\goth g} \mapsto \hat{X} \in Vect(P). 
\]
Hence the horizontal subbundle ${\mathcal H}P$ can be described not only
as a kernel of the projection operator
 $Pr: {\mathcal T}\!P \to {\mathcal V}\!P$, but also 
as a kernel of a  connection 1--form:
\[
\lambda: {\mathcal T}\!P \to {\goth g}.
\]
 It is a $G$--invariant form on the principal $G$--bundle $P$ with values in
the Lie algebra $\goth g$, such that
 \[
 \imath_{\hat{X}} \lambda = X, \quad
\text{for } \quad X \in {\goth g}. \]
%%%%%%%%%%%%%%%%%%%%%%%%%%%%%%%%%%%%%%%%%%%%%%%%%%%%%%%%%%%%%%

 Now let $G$  act on a symplectic
manifold $(F, \sigma)$ by symplectomorphisms, i.e. there is a 
group homomorphism
\[
G \to Symp (F,\sigma): \quad g \mapsto \psi_g.
\]
The infinitesimal action determines the Lie algebra homomorphism 
\begin{equation}
 \label{eq:vf}
{\goth {g}} \to Vect (F, \sigma): X \mapsto \hat{X}, \  \ \ 
\text{defined by} \  \ \
\left. \hat{X} = \frac {d}{dt}\right| _{t=0}
 \psi_{exp (t X)} 
\end{equation}
 A symplectic form determines a correspondence between functions
and certain vector fields, called Hamiltonian vector fields:
 \begin{equation}
 \label{eq:f-vf}
{\C^\infty(F)} \to Vect (F, \sigma):H \mapsto X^H, \  \ \ 
\text{defined by} \  \ \
\imath_{{X}^H} \sigma = d H.
\end{equation}

\begin{definition} 
\label{def-ham}
The action of $G$ on $F$ is called Hamiltonian if:
\begin{enumerate}
\item
There is a lift ${\goth g} \to 
 C^\infty (F): \quad  X \mapsto H_X$
 \[
   \dgARROWLENGTH=0.6\dgARROWLENGTH
   \begin{diagram}                 
 \node[2]{C^\infty(F)}\arrow{s}{}	\\
\node{\goth g} \arrow{ne,..} \arrow{e}{}      \node{Vect(F,\sigma)}
   \end{diagram}
\]
  This means that
 there is a Hamiltonian function
$H_{X}$ so that 
 $ \imath_{\hat{X}} \sigma = d H_X.$
 \item 
This map is a  Lie algebra homomorphism:
\[ 
H_{[X,Y]} = \{H_X, H_Y\}.
\]  
$[X,Y]$ is the Lie bracket  in $\goth{g}$ and 
$ \{H_X, H_Y\}$ is  the Poisson bracket in $C^\infty (F)$.
\end{enumerate}
\end{definition}
\noindent (If a group action satisfies only first condition it is called
weakly Hamiltonian.)
 Let us also mention the following equality:
 \begin{equation}
\label{eq:xy}
H_{[X,Y]} = \hat{X} H_Y - \hat{Y} H_X
 \end{equation}
 Hamiltonian action determines a map
$\mu: F \to {\goth g}^*,$
for each point $x \in F$ defined by 
 \[
<\mu (x)| X> = H_X (x)
 \]
 for any
 $X \in {\goth g}$, where $< \cdot | \cdot >$ is the pairing: ${\goth
g}^* \times {\goth g} \to \C$. This map is usually called a moment map,
however for our purposes of quantization we will call a map ${\goth g}
\to C^\infty (F)$ also a moment map and it is this map of algebras which
we are going to quantize.

\noindent {\bf Weak coupling: connection $\longleftrightarrow$
symplectic form.}  The following proposition in its present form is an
adaptation for our purposes of a theorem about weak coupling form from
\cite{GLS}. 
 \begin{proposition}
	\label{prop:form}
Let $G\to Symp (F,\sigma): g \mapsto \psi_g$ be a Hamiltonian
action on $(F,\sigma)$ with a moment map $\mu^F$. 
Then every connection $\Gamma$ on the principal $G$--bundle 
$P \to B$ over a symplectic manifold  $(B, \omega^B),$ 
gives rise to a one--parameter 
family of symplectic forms on the associated bundle
$M = P\times_G F \to B$, which restricts to the forms 
$\sigma_b$ on the fibres:
\begin{equation}
\Omega_{\epsilon} =  {\epsilon^2} \omega^\Gamma +  \pi ^* \omega ^B
   \label{eq:form}
\end{equation}
where  $\epsilon$ is a small parameter and $\omega^\Gamma $
is the coupling form, so that at a point $x \in M$, $\pi (x) = b$
\[  
\omega^\Gamma = \sigma_b + H_T 
\]
 where $T \in {\cal A}^2 (B)$: $T(X,Y) = -Pr([X^H,Y^H]), \ X,Y \in
{\mathcal T}B $ and $H_V $ is a Hamiltonian function of a vector field
$V$ with respect to the form $\sigma_b$, defined by $\mu_F$.
 \end{proposition}

\begin{remark}
 Notice that $\sigma_b$ is nonzero only on vertical vectors, while $H_T$
is nonzero only on horizontal ones.  This extra term $H_T$ is needed for
the form to be closed, $\epsilon $ makes the form $\Omega_{\epsilon} $
to be nondegenerate.
 \end{remark}

%%%%%%%%%%%%%%%%% 
\begin{pf}  
Sketch (for the full proof see \cite{GLS} or \cite{McDS}).
The main idea is to use the so called 
Weinstein universal phase space --
\[
W = P \times {\goth g}^*.
\]
 Given a connection on $P$ 
the space $W$ could be identified with 
the  vertical subbundle of the cotangent bundle
 \[
W = P\times_G {\mathcal T}^*\!G = {\mathcal V}^*\!P.
\]
 A connection is  a splitting 
$ \Gamma: {\mathcal T}\!P  =  {\mathcal H}\!P \,  \oplus \,  {\mathcal V}\!P $ 
and ${\mathcal V}^*\!P$ is defined as one--forms 
which vanish on horizontal vectors:
${\mathcal V}^*\!P = ( {\mathcal H}\!P)^\perp.$
 Hence it has a $G$-equivariant 
symplectic form coming from 
the canonical symplectic form on  
${\mathcal T}^*\!P$. Moreover, the 
action of the group $G$ on $W$ is Hamiltonian (see for example \cite{Ar}).

 The moment map
$
\mu^W: W \to  {\goth g}^*$
is given by the projection
\[
pr_{{\goth g}^*}: W \to {\goth g}^*.
\]
Then the symplectic reduction of $W \times F$ at $0$ value of the moment map
$\mu = \mu^W +\mu^F$
is exactly $M = P\times_G F$, and the symplectic form on $M$ is 
inherited from $W$.
  
The explicit formula is obtained in the following way.
 Let the connection $\Gamma$ be given by a connection $1$--form,
$\lambda_p :{\mathcal T}\!_p P \to {\goth g}.$  It determines
a horizontal subbundle in ${\mathcal T}P$ by $ 
{\mathcal H}_p P = \{ v \in {\mathcal T}\!_pP \mid \lambda_p (v) =0 \}.$

\noindent ${\mathcal V}^*\!P = ({\mathcal H}P)^{\perp}$ 
is also defined by $\lambda$.
The connection $\lambda :{\mathcal T}\!P \to {\goth g} $ and together with 
the action
$\rho: {\goth g} \to {\mathcal V}\!P$ define the 
injection 
 \[
\imath_\lambda: {\mathcal V}^*\!P \hookrightarrow {\mathcal T}^*\!P. 
 \]
 By definition of the connection $1$--form, this injection is 
equivariant under the action of $G$ and hence the $2$--form 
\[
\omega_\lambda = \imath_\lambda^*  \ \omega_{can} 
\in {\mathcal A}^2({\mathcal V}^*\!P)
\]
is invariant under the action of $G$.
 This pull--back of the canonical symplectic form on 
${\mathcal T}^*\!P$ gives
a closed $2$--form on ${\mathcal V}^*\!P$. 

\noindent The canonical $1$--form $\alpha$ on ${\mathcal T}^*\!P$ 
is defined as follows. Let $(p,s_p)$ be a point in ${\mathcal T}^*\!_pP$,
let also $v$ be a tangent vector field  in the tangent bundle 
$\pi: 
{\mathcal T}\!({\mathcal T}^*\!P) \to {\mathcal T}^*\!P$ then at the point 
$(p,s_p)$
\[
<\alpha |v >_{(p,s_p)} = - <s_p |\pi_* v>_p.
\] 
Then the pullback of the canonical one--form to ${\mathcal V}^*\!P$ is 
$<\imath_\lambda^* \alpha| v>_{(p,s_p)} = - < pr_{{\goth g}^*}s_p |
 \lambda (\pi_* v) >_{p}$ so 
\begin{equation}
\label{eq:prform}
\omega_\lambda = - d < pr_{{\goth g}^*} | \lambda >.
\end{equation}
 Now the form on $W \times F$ is $\omega_\lambda + \sigma$. 
We now want  
the form on the reduced space  $M = (W \times F )// G$ 
at the regular value of the 
moment map $\mu^W + \mu^F = 0$. Since  $ \mu^W =  pr_{{\goth g}^*}$
 the form on $M$ becomes  
\[
\omega^\Gamma = d < \mu^F |
 \lambda > + \sigma.
\]
 The one--form $< \mu^F | \lambda >$ can be rewritten as the
Hamiltonian $ H_\lambda $. It should be understood  in
the following way. The connection $1$--form $\lambda: {\mathcal T}\!P
\to {\goth g}$ defines the connection on the associated bundle $M = P
\times_G F$. The horizontal subbundle in ${\mathcal T}\!M$ is 
the image of ${\mathcal H}\!P$ under the map $P \times F \to M$. By
abuse of notation we call the map ${\mathcal T}\!M \to {\goth g}$ also
$\lambda$. From now on $\lambda$ is a $1$--form on $M$ with values in
the Lie algebra ${\goth g}.$ Hence $ H_\lambda $ is a $1$--form on $M.$

Its differential gives a two form 
 $d < \mu^F | \lambda >= d (H_\lambda)$. Applied to two vectors 
$V, W \in {\mathcal T}\! M:$
\[
 d (H_\lambda) (V, W)= V H_\lambda (W) - V H_\lambda (W) - H_\lambda([V,W])
\]
using (\ref{eq:xy}) we get that it is nonzero only on horizontal vectors and 
gives  the Hamiltonian of the curvature.
We see also that $\omega^\Gamma$ restricted to fibres gives the
symplectic form on the fibres.
 \end{pf}

This construction is quite general \cite{GLS}:
The symplectic fibrations with connection constructed this way turn out
to include all symplectic fibre bundles with connection for which 
the holonomy group is a finite dimensional Lie group. 

%%%%%%%%%%%%%%%%%%%%%%%%%%%%%%%%%
\noindent {\bf In local coordinates.} 
 \label{loc}
Let us take a point $x \in M$. One can introduce a local frame $\{f_\alpha \}$
of vertical tangent bundle ${\mathcal V}M$ and a local frame $\{e_i\}$ in 
${\mathcal T}B$ at
a point $b  = \pi( x)$ of $B$, with dual frames $\{f^{\alpha}\}$ and 
$\{e^i\}$. Using the connection we obtain a local frame on the
tangent bundle ${\mathcal T}M = \pi^*{\mathcal T}B \oplus {\mathcal V}M$ 
at a point $x$. 

\noindent
Then the form can be written as a block matrix:
\begin{equation}
\Omega_\epsilon= 
\left| 
	\begin{array}{cc}
 {\pi ^* \omega ^B } + {\epsilon^2} H_T & 0 \\
              0                        & {\epsilon^2}\sigma_b
        \end{array}
\right|
\label{eq:formM}
\end{equation}
Hence the corresponding Poisson bracket is also a block matrix:
\[
\left| 
	\begin{array}{cc}
 ({\pi ^* \omega ^B } + {\epsilon^2} H_T)^{-1} & 0 \\
             0                         & {\epsilon^{-2}} \sigma_b^{-1}
        \end{array}
\right|
\]
We see that the Moyal product with respect to this form is 
a product of those on  the base and on the fibres. 

As for the connection, let 
$\Gamma^{\delta}_{\beta \gamma}$ be the Christoffel symbols of the 
symplectic connection on ${\mathcal T}F_b$ preserving the form $\sigma_b$ 
along the fibre $F_b$. This connection can be written in the 
	²local coordinates as follows:
\begin{equation}
\nabla^F = d^F +  \disph 
[\Gamma_{\alpha \beta \gamma} f^\alpha f^\beta d \xi^\gamma, \cdot ].
\label{eq:nablaF}
\end{equation}
 where  $\Gamma_{\alpha \beta \gamma} 
= (\sigma_b)_{\alpha \delta} \Gamma^{\delta}_{\beta \gamma}$. 
%%%%%%%%%%%%%%%%%%%%%%%%%%%%%%%%%%%%%

\section{ Quantum moment map.\\}
\label{qmm}
 Let $(F, \Sigma,  \nabla)$ be an $F$--manifold.
Here $\Sigma$ is a deformation of the symplectic form  $\sigma$:
 \[
 \Sigma = \sigma + \hbar \sigma_1 + \hbar^2 \sigma_2 + \cdots,
 \]
 $\sigma_i$ being  closed 
$2$--forms on $F$.
Let ${\A} (F)$  be the corresponding quantization of $F$ 
with the characteristic class $[\Sigma] \in
H^2(F)[[\hbar]]$.
 ${\A} (F)$  is a  noncommutative algebra 
of formal series in $\hbar$ with coefficients being 
smooth functions on $F$. 
 The $\ast$--product on ${\A}(F)$ defines the Lie algebra structure: 
 \[
[f,g]_\ast = \disph (f \ast g - g\ast f), 
\quad \text{for} \ \ f,g \in {\A} (F).
\]
 Let $G$ be a group acting on  
$F$ Hamiltonially (Definition \ref{def-ham}).
Let ${\goth g} \to C^\infty (F)$ be its  moment map.
 There is an induced action 
of $G$ on  ${\A} (F)$.  
We want to quantize the moment map,
namely, get a Lie algebra map from the algebra ${\goth g}$ to the
quantized algebra ${\A} (F)$. However it is possible only up-to a
two--cocycle in $C^\infty (F)[[\hbar]]$, so we get a projective
representation or  instead we 
should consider a central extension of ${\A} (F).$ We could also 
slightly change a definition of a quantum moment map. 
We eliminate the 
central elements by considering the map into the adjoint representation 
of ${\A} (F)$, the inner automorphisms $Inn \ {\A}(F)$. They obviously 
inherit the Lie algebra structure from ${\A} (F)$, so 
 \begin{definition}
 A quantum moment map is a map of Lie algebras 
 \[
\mu^{\vee}: {\goth g} \to Inn \ {\A}(F)
\]
 with the correspondence principle:
$
\mathop{lim_{\hbar \to 0}} \mu^{\vee}(X)(f) = \{ H_X, f \}
$
for $X \in {\goth g}$ and $f \in {\A}(F)$.
 \end{definition}
\begin{remark}
This definition could be 
reformulated through an isomorphism of associative algebras:
$
 \mu: U({\goth g}) \to {\A}(F),
$
such that on vector fields it gives:
\[
 \mu(X)^{\vee}(f) = [\mu(X), f]_\ast.
\]
 \end{remark}
 \begin{proposition}
 \label{prop:action}
 Let $G$ act Hamiltonially on
$(F, \sigma)$ with the Hamiltonian function $H_X$.  Let 
\[
(F, \Sigma=  \sigma + \Sigma_{i=1}^{\infty}\hbar^i \sigma_1 , \nabla)
\]
be an $F$--manifold such that $\nabla$ is a $G$--invariant connection
and assume that one can define functions $H^i_X, i= 1,2, \cdots$  as follows 
\begin{equation}
\label{eq:Hi}
\iota_{\hat{X}} \sigma_i = d H^i_X.
\end{equation}
Let $ {\A} (F)$
be the algebra of quantized functions. 
Then the quantum moment map is given by 
 \begin{equation}
 \label{eq:qmoment}
 \mu(X) = H_X + \Sigma_{i=1}^{\infty} \hbar^i {{H}^i}_X
\end{equation}
Also ${\mu}^{\vee}: {\goth g} \to Inn \ {\A}(F), \ 
\mu(X)^{\vee}(f) = [\mu(X), f]_\ast$ 
is a homomorphism of Lie algebras:
 \begin{equation}
 \label{eq:qbracket} 
{\mu}^{\vee} ([X,Y])= [{\mu}^{\vee} (X), {\mu}^{\vee} (Y)]_\ast.
\end{equation}
Moreover, there are no higher terms in $\hbar:$
 \begin{equation}
\label{eq:mu}
 [\mu(X), f]_\ast = [H_X, f]_\ast = \{ H_X, f \}.
 \end{equation}
\end{proposition}
 \begin{pf}\footnote{I am grateful to Boris Tsygan for pointing out 
a gap in the proof in an earlier version of the paper.}
 We are going to prove (\ref{eq:qmoment}--\ref{eq:mu}) by lifting the
Hamiltonian $H_X$ to a section of the Weyl algebra bundle, since the
$\ast$--product on $C^\infty (F)$ is defined through the the Weyl
algebra bundle (\ref{eq:Q}) and (\ref{eq:ast}). Namely, let
$\Gamma_{\text{D--flat}} (F, \W_F)$ be the space of flat sections 
Fedosov connection $D$ constructed from $\nabla$ corresponding to the
quantization of $F$ with the characteristic class $\Sigma$. Then there
is a one--to--one correspondence
 \[
Q: C^\infty (F)[[\hbar]] \to \Gamma _{\text{D--flat}}(F, \W_F)
\]
 which defines a product on $F$ 
 \[
f\ast g = Q^{-1} (Q(f) \circ Q(g)).
\]
 The structure of a Lie algebra in 
 Weyl algebra bundle ${ \W_F}$ is defined by the  fibrewise 
commutator 
\[
[{\bold a} ,{\bold b }]_\circ ={\bold a} \circ {\bold b} - {
\bold b}
\circ {\bold a}, \quad \mbox{for} \ {\bold a},{\bold b } 
\in \Gamma (F, \W_F)
\]

  Recall  that  the  map
 $H: {\goth{g}} \to C^\infty (F)$ is given by the condition 
\begin{equation}
\iota_{\hat{X}} \sigma = d H_X, 
 \label{eq:act}
 \end{equation}
 where $\hat{X}$ is a field corresponding
to $X$ under the map ${\goth g} \to Vect(F) $ (\ref{eq:vf}).
 We find the image of a Hamiltonian 
 $H_X$ in $\Gamma _{\text{D--flat}}(F, \W_F)$
 generalizing to the case of a deformed symplectic form the proof of  
Fedosov \cite[Propositions
5.8.1,2]{F2},\cite{F3}. 
By construction the Fedosov connection
$D$ is $G$--equivariant since any element of
the group $G$ preserves the initial symplectic connection $\nabla$.  An
easy calculation shows also that flat sections of an equivariant
connection $D$ are also equivariant.
  It also means that the Lie derivative of $D$ is zero with respect to
any vector field $\hat{X}:$
 \[ 
[L_{\hat{X}}, D]_\circ = 0.
\]
 Since $L_{\hat{X}}$ and
$D_{\hat{X}}= \imath_{\hat{X}} D + D \imath_{\hat{X}}$ are first order
derivations commuting  with $D$ we can find an analogue of the Cartan
homotopy formula for the Lie derivative on forms with values in the Weyl
algebra bundle. The difference of $L_{\hat{X}}$ and $D_{\hat{X}}$
could  only be an inner automorphism of the Weyl algebra bundle, 
we denote it  $[ Q(X), \cdot ]_\circ:$
 \begin{equation}
\label{eq:Lie}
L_{\hat{X}} {\bold a} = (\imath_X D +  D  \imath_X) {\bold a} 
+ [ Q(X),{\bold a}]_\circ.
\end{equation}
It is easy to see that 
the equality  $(\ref{eq:Lie})$ is true in Darboux coordinates chosen 
so that the field $\hat{X}$ is just a pure derivation 
in the direction of only one of the coordinates.
Then $ D = D^0 =  d + \delta$ and 
 $Q(X)= Q^0(X) = \mbox{\it central section}(X, \hbar) - \imath_X \delta.$
 Let $D = D^0 + [\Delta \gamma, \cdot]_\circ$ be another flat
connection, $\Delta \gamma$ being an equivariant one form in ${\cal
W}_F$.  Then since the $L_{\hat{X}}$ does not change and commutes with
the new $D$ as well, the right hand side of (\ref{eq:Lie}) must not
change either, so one has to subtract $\imath_X \Delta \gamma$ from
$Q^0(X).$

We want to show that we could chose a  $\mbox{\it central section}(X,
\hbar)$ in $Q(X)$  to be equal to $\mu (X)= H_X  + \hbar \cdots$ so  that 
$Q(X) =\mu (X) - \imath_X \gamma$
becomes a quantization of the moment map, that is a flat section
corresponding to $H_X.$

 Recall that locally we can write a  connection
$D$ as $D f = d f + \disph [\gamma, f]_\circ$ (\ref{eq:con}). 
The equation $D^2 = 0$ becomes (\ref{eq:eq}):
 \begin{equation}
\label{eq:sigma}
\Sigma + d \gamma + \disph \frac{ [\gamma, \gamma ]_\circ}{2} = 0.
 \end{equation}
Then  since $L_{\hat{X}} D =0$ from (\ref{eq:Lie}) 
\[
0 = \{ \imath_{\hat{X}} (d + \disph  [\gamma, \cdot]_\circ) + 
(d +\disph [\gamma, \cdot]_\circ) \imath_{\hat{X}}  \} (D) 
- [ \imath_X \gamma, D]_\circ
= (\imath_{\hat{X}} d + d \imath_{\hat{X}}) (D)
\] 
 we get
\begin{equation}
\label{eq:L}
(\imath_{\hat{X}} d + d \imath_{\hat{X}} ) \gamma = 0.
\end{equation}

 Using (\ref{eq:Hi}),(\ref{eq:act}), and also  
(\ref{eq:sigma})  and (\ref{eq:L}) we find  that
indeed $Q(X)$ is a flat section of $D:$
\begin{eqnarray*}
D Q(X) &= 
d \{ H_X + \hbar H^1_X + \cdots \}-  d(\iota_X \gamma)
+  \disph [\gamma,  \{ H_X + \hbar H^1_X +  \cdots \}]_\circ
-  \disph [\gamma,\iota_X \gamma]_\circ\\
& = \imath_X \Sigma  + \imath_X d \gamma+ 0 +\imath _X 
(\disph \frac{ [\gamma, \gamma ]_\circ}{2}) 
 = \imath_X (\Sigma +  d \gamma + \disph \frac{ [\gamma, \gamma ]_\circ}{2}) = 0.
\end{eqnarray*}
 We also get that $\mu(X) = H_X + \Sigma_{i=1}^{\infty} \hbar^i
{{H}^i}_X$, so the commutator
 \begin{eqnarray*}
[\mu(X), f]_\ast & = \disph Q^{-1}\{  [Q(X), Q (f)]_\circ \}
= \disph Q^{-1} \{ [  (H_X + \hbar H^1_X + \cdots - \imath_{\hat{X}} \gamma),
Q (f) ]_\circ \}\\
 &=\disph Q^{-1}\{  [ - \imath_{\hat{X}} \gamma,
Q (f) ]_\circ \} = Q^{-1} \{ - \imath_{\hat{X}} (D - d)Q (f) \}
=  Q^{-1}\{ \imath_{\hat{X}}  d Q (f) \} \\
& = L_{\hat{X}} f = {\hat{X}} (f) = \{ H_X, f \}.
\end{eqnarray*}
Then since the action is Hamiltonian we get  on the
quantum level:
 \begin{equation}
 \label{eq:mu2}
[\mu^{\vee}(X),\mu^{\vee}(Y)]_\ast = ad [H_X, H_Y]_\ast = 
ad (\{H_X, H_Y\}) = ad (H_{[ X, Y]}) = \mu^{\vee}([ X, Y]).
 \end{equation} \end{pf}

\noindent For the case $\Sigma = \sigma$, in the absence of additional terms 
$\hbar^i H^i$ we do not need to consider the adjoint representation.  
The notion of the quantum moment map repeats the classical one,
it is a map from ${\goth{g}}$ to ${\A}(F):$
\begin{corollary}
Consider a quantization $ {\A}(F)$ obtained from the $F$--manifold
$(F, \sigma, \nabla)$. Let $G$ act Hamiltonially on
$F$ with the Hamiltonian function $H_X$.
 Then the quantum moment map is a homomorphism of Lie algebras
${\mu}: {\goth g} \to ({\A}(F)):$
\[
 \mu(X) = H_X 
\quad \text{and} \quad
[\mu(X), f]_\ast = \{ H_X, f \}.
 \]
\end{corollary}

\begin{remark} We got that 
  $\mu^\vee:{ \goth {g}} \to Inn \ {\A}(F)$ is a map of Lie
algebras. This gives the  positive answer to a question posed in \cite{Xu}
that every classical moment map can be uniquely lifted to a quantum
moment map to ${\A}(F)$. 
\end{remark} 

\section{Quantization of twisted products. \\}
\label{chap:main}

\subsection{Auxiliary bundle: quantization of the fibres.\\}
\label{sect:aux}
In this section we look at the bundle of quantizations along the fibres
of our symplectic fibration $ M \to B$ and construct a covariant derivative 
on this quantized bundle from the connection on the  bundle $ M \to B$.

\noindent {\bf Fibrewise quantization.}
Again, let $\pi: M \to B$ be a locally trivial fibration with a fibre
$F$. Let $F$ be equipped with a $F$--manifold $(F, \Sigma,  \nabla)$.
Here we take $\Sigma$ to be arbitrary characteristic class, i.e. 
cohomology class $[\Sigma] \in H^2(F)[[\hbar]]$.
 \[
\Sigma = \sigma + \hbar \sigma_1 + \hbar^2 \sigma_2 + \cdots,
\]
where $\sigma$ is a symplectic form on $F$.

Let ${\A} (F)$ be the quantization of $F$, i.e. noncommutative algebra 
of formal series in $\hbar$ with coefficients being 
$C^\infty$--functions on $F$.
We are  defining the  bundle  ${\U}$ over $B$ such that its 
  fibres  are algebras
of quantized functions on the fibres of the bundle $\pi: M \to B$, i.e.
the fibre of ${\U}$ at a point $b \in B$ is
\[
{\U}_b = {\A} (M_b).
\]
 The structure group $G$ of  $\pi: M \to B$
acts on $F$ by symplectomorphisms hence there is  a $G$--action on ${\A} (F)$.
Since $M = P \times_G F$ is  an associated bundle to a $G$--bundle $P \to B$, 
this auxiliary  bundle $\U$ is  also associated to 
$P$ with the fibre ${\A} (F)$
\[
{\U} = P \times_G {\A} (F).
\]

\noindent {\bf Covariant derivative on the auxiliary bundle.}
\label{sec:derivative}
 If a connection on a fibre bundle is defined, 
there is a geometric way of putting a connection on the infinite
dimensional bundle of smooth sections over fibres of a given
finite dimensional fibre bundle.

 In our case, covariant derivative on the bundle ${\U} \to B$ respecting
the algebra structure can be obtained from a connection 1--form on the
principal bundle. Understanding of the formula for this covariant
derivative is important for the sequel.  We also find that the curvature
of this covariant derivative equals to the 
adjoint action of  the first summand
in the coupling form.
  
In general, a choice of a connection on the principal $G$--bundle $P$
determines a covariant derivative on any associated vector bundle (we
follow an exposition in \cite{BGV}).

Let $\lambda \in {\mathcal A}^1(P, \goth g)$ be a connection 
1--form on a  principal $G$--bundle  $P$.
Let also $G$ act on some vector space $E$, 
and the action be given by the map $ \rho,$
\[
\rho: G \to End( E).
\]
Then the bundle   ${\cal E} = P\times_G E \to B$ 
is an associated bundle to 
the principal bundle $P$. The space of differential 
forms on $B$ with values in ${\cal E}$,
${\cal A}^k (B,{P \times_G E})$, can be described as the  subspace of the 
space of differential forms on $P$ with values in $E$.
This subspace is a space of all basic forms with values in $E$, 
${\cal A} (P,E)_{bas}$.
A {\em basic differential form} on a principal bundle $P$ with a 
structure group $G$, taking values in the representation $(E,\rho)$ of $G$,
is an invariant and horizontal differential form, that is a form 
$\alpha \in {\cal A} (P,E)$ which satisfies
\begin{enumerate}
\item $g \cdot \alpha = \alpha, \quad g\in G$ 
\item
$\imath (X) \alpha = 0 $ for any vertical field $X$ on $P$.
\end{enumerate}
\begin{lemma} 
If $\alpha \in {\cal A}^q (P,E)_{bas}$, define 
$\alpha_B \in {\cal A}^q (B, P \times_G E)$ by
\[
\alpha_B (\pi_*X_1, \ldots , \pi_*X_q)(b) = 
[ p, \alpha(X_1, \ldots ,X_q) (p)],
\]
where $p \in  P$ is any point such that $\pi(p) =b$, 
and $X_i \in {\mathcal T}_pP$.
Then $\alpha_B$ is well defined, and the map 
$\alpha \to \alpha_B$ is an isomorphism from
${\cal A}^q (P,E)_{bas}$ to ${\cal A}^q (B, P \times_G E).$
\end{lemma}

As a particular case, there is a representation of the sections of
$\cal E$ as $G$--equivariant functions on $P$ with values in $E$.
 Let ${C}^\infty (P,E)^G$ denote the space of equivariant 
maps from $P$ to $E$, 
that is those maps $s: P \rightarrow E$ that satisfy 
$s(p \cdot g) = \rho(g)s(p)$. There is a natural isomorphism between
$\Gamma (B, P\times_GE)$ and ${C}^\infty (P,E)^G$, given by sending 
$s \in {C}^\infty (P,E)^G$ to $s_B$ defined by 
\[
s_B(b)= [p, s(p)];
\]
here $p$ is any element  of $\pi^{-1}(b)$ and $ [p, s(p)]$ is the element of 
${\cal E} = P \times_G E$ corresponding to  $(p, s(p))\in {C}^\infty (P,E)^G$. 
Infinitesimally, a function $s$ in  ${C}^\infty (P,E)^G$ satisfies 
the formula:
\[
(X_P \cdot s)(b) + \rho (X) s(b) = 0, \ \ \text{for} \ X \in {\goth g},
\]
where we also denote by $\rho$ the differential of the representation $\rho$:
\begin{equation}
 \label{eq:rho}
          \rho: {\goth g} \to Vect( E).
\end{equation}
Given a connection $1$--form $\lambda$ on $P$ one obtains the
covariant derivative $\nabla$ on the  associated 
 vector bundle $\cal E$ from the following commutative diagram: 
\begin{equation}
\begin{diagram}
\node{{C}^\infty (P,E)^G } \arrow[2]{e,t}{d+ \rho (\lambda)}
      \arrow{s,l}{} 
        \node[2]{{\cal A}^1 (P,E)_{bas}} \arrow{s,r}{}\\
\node{\Gamma (B,\cal E)} \arrow[2]{e,t}{\nabla}
   \node[2]{{\cal A}^1 (B,\cal E)}
		\label{eq:diag2}
\end{diagram}
\end{equation}
 In other words, the covariant derivative can be written 
$(\nabla_X s) (b)= [p, (X^H s_P)(p)]$, where $X \in  {\mathcal T}B$
and $X^H$ is its horizontal lift to the principal bundle $P$. 
The formula for the covariant derivative 
on our auxiliary bundle  bellow should be understood  by means of the 
diagram (\ref{eq:diag2}).

 Now let us return to our particular case, namely:
\[
{\cal E} = {\U} = P \times_G {\A} (F), \ \text{so } \  E = {\A} (F),
\]
 Then (\ref{eq:rho}) becomes a map $\rho: {\goth g} \to Vect ({\A} (F))$,
or in other words it is given by the moment map $\mu^{\vee}: {\goth g}
\to Inn \ {\A} (F)$.
 \begin{proposition}
\label{th:con} 
Covariant derivative on 
${\U} \to B$ 
corresponding to a connection 
1--form $\lambda$ on $P \to B$ is given by the formula:
\begin{equation}
\nabla^{\F} f = df + [H_{\lambda}, f]_\ast = df +\{H_{\lambda}, f\}.
						\label{eq:auxcon}
 \end{equation}
Its curvature is a $2$--form on $B$
 with values in  ${\U}$:
\[
R^{\F}f =  \{H_T, f\}, \quad \text{where} \ \
 T (X,Y) = - Pr {[{X}^H,{Y}^H]}, 
\quad Pr : {\mathcal T}P \to {\goth g}.
\]
\end{proposition}
\begin{pf}
The covariant derivative formula  
follows from  (\ref{eq:act})  and the diagram (\ref{eq:diag2}). 
 Then the general definition of 
the curvature of  a covariant derivative 
$\nabla: \Gamma(B, E) \to {\mathcal A}^1(B, E)$ is
 \[
R (X,Y) = \nabla_{\tilde{X}} \nabla_{\tilde{Y}} - 
\nabla_{\tilde{Y}} \nabla_{\tilde{X}} - 
\nabla_{[\tilde{X},\tilde{Y}]}, 
\]
where $ X, Y \in {\mathcal T}B,$ and $ \tilde{X}, \tilde{Y}$ 
are their horizontal lifts. In our particular case the 
expression for the curvature follows from the formulas (\ref{eq:mu}) and 
(\ref{eq:mu2}). 
\end{pf}

%%%%%%%%%%%%%%%%%%%%%%%%%%%%%%%%%%%%%%%%%%%%%%%%
\subsection{Fedosov connection and flat sections on symplectic fibrations.\\}
 In this section we will show that the complex of differential forms
with values in some twisted Weyl bundle gives a resolution of
$\Gamma(B, \U)$.

\noindent {\bf Bundle of sections.}  The structure of the
bundle $\pi:M \to B$ is reflected in the representation of the space of
functions on $M$ as sections of a certain bundle over $B$.  Namely, let
${\cal F} $ be the bundle over $B$, such that a fibre over $b \in B$ is
the space of series in $\hbar$ with coefficients being functions on the
fibre $M_b$ of the bundle $\pi:M \to B$:
 \[
{\cal F}_b = C^\infty(M_b)[[\hbar]]. 
\]
Then 
$C^\infty(M) [[ \hbar]]$ can be represented as sections of 
the bundle ${\cal F} $:
\begin{equation}
\label{eq:C(M)}
C^\infty(M) [[ \hbar]] = \Gamma(B, {\cal F})
\end{equation}
Hence, we can try to obtain a  quantization of $M$ by 
quantization of  the bundle ${\cal F} \to B$.
 This leads us to consider the twisted Weyl algebra bundle over $B$.

Indeed, let ${\cal W}_M$  be  the Weyl algebra bundle on $M$. Then consider 
${\cal W}_{M/B}$, the  bundle on $B$ with the fibers ${\cal W}_b$
being the restriction of the bundle ${\cal W}_M$ to
the fibre  $M_b$, $b \in B$. Given a 
symplectic connection (\ref{eq:split}) on the bundle 
$\pi: M \to B$ we can talk about the following 
 isomorphism on the level of sections 
 \begin{equation}
\beta: \Gamma(M,{\cal W}_M)  \to 
\Gamma(B,{\cal W}_B \otimes_{\C [[\hbar]]} \Gamma( M/B,{\cal W}_{M/B}))
\label{eq:M/B}
	\end{equation} 
or more generally for differential forms:
\begin{lemma}
Given a symplectic connection:  
\[
 \quad {\mathcal T}\!M = {\mathcal H}M 
\oplus {\mathcal V}\!M, \quad {\mathcal H}_zM \cong
 \pi^* {\mathcal T}\!_{\pi(z)}B
\]
there is an isomorphism:
 \begin{equation}
	{\cal A}^n (M,{\cal W}_M) \to \oplus_{p+q=n}{\cal A}^p (B, {\cal W}_B
\otimes _{C [[\hbar]]} {\cal A}^q ( M/B, {\cal W}_{M/B}))
				  \label{eq:M/B1}
	\end{equation} 
\end{lemma}
\begin{pf} At each point $z$ of $M$ the Weyl algebra can be defined
as the universal enveloping algebra of the Heisenberg algebra 
of ${\mathcal T}^*_zM \oplus \hbar {\R}$. Universal enveloping algebra 
is by definition a quotient of the tensor algebra by a certain ideal.
The ideal is 
\[
I= \{ e \otimes f - e \otimes f  + i \hbar \omega^{-1} (e,f) \}, \quad \text{for}
\quad e,f \in {\mathcal T}^*_zM
\]
 and $\omega$ the symplectic form on  ${\mathcal T}^*_zM $.
The connection on ${\mathcal T}M $ splits the ideal into a  sum
of a horizontal and vertical ideals:
\[
I = I_H + I_V.
\]
 Thanks  to the splitting of the symplectic form (\ref{eq:formM})
vertical ideal $I_V$ is also an ideal in  the 
tensor algebra of ${\mathcal V}^*\!M$, which leads to the result.
\end{pf}
\begin{remark}
\label{r:iso} Notice that we were talking about the isomorphism of
bundles, but we did not show yet that it is an isomorphism of algebras 
-- it is the subject of the section \ref{sec:qts}.
\end{remark} 
 We see that one can consider fiberwise quantization as a first step in
quantization of the total space. This leads to the auxiliary bundle
described  in Section \ref{sect:aux}. 
The quantization map along the fibres:
\begin{equation}
 \label{eq:Qf}
Q^{\text{fibre}}: \Gamma(B, {\cal F}) \to \Gamma(B, {\U})
\end{equation}

%%%%%%%%%%%%%%%%%%%%%%%%%%%%%%%%%%%

 \noindent {\bf Twisted Weyl algebra bundle.}
Let $ \W_B$ be the Weyl algebra bundle corresponding to the $F$--manifold
 $\{B, \omega, \nabla^B \}$. Consider the twisted 
bundle  
\[
{\B} = \W_B \otimes_{\C[[\hbar]]} {\U} \to B.
\]
Thanks to the splitting of the form (\ref{eq:formM}) 
the  sections of  $\W_B$ and sections of ${\U}$ 
commute with each other.

\begin{remark} \label{rem} $\Gamma(B,{\B})$ can be considered as a space
of fiberwise flat sections of the bundle obtained by $\beta$
(\ref{eq:M/B}): $\Gamma (B,{\cal W}_B \otimes_{\C [[\hbar]]}
 \Gamma(M/B,{\cal W}_{M/B}))$. Indeed, the $F$--manifold $(F, \sigma,
\nabla^F)$ induces a structure of $F$--manifold on each fibre. It
defines corresponding isomorphism of flat sections of Weyl algebra
bundles:
 \begin{equation}
\label{eq:U}
 \U  = \Gamma_{\text{flat}}( M/B,{\cal W}_{M/B})
 \end{equation}
\end{remark}

This bundle ${\B} $ is a bundle of graded algebras with degrees assigned
as in the original $\W_B$ bundle, namely,
 \begin{equation}
deg (\hbar) = 2,  \quad 
deg (e^i)=1, \quad \text{$e^i$'s being generators of $\W_B$}.
\label{eq:deg}
\end{equation}
Let 
\[
({\B})_n =  \{ s \in {\B}, \quad \text{such that} \quad 
deg(s) \geq n \}.
\]
Then let us also define the grading
\[
gr_n ({\B})= \{ s \in {\B}, \quad \text{such that} \quad 
deg(s) = n \}. 
\]
It is isomorphic to $({\B})_n / ({\B})_{n+1}.$
 The pointwise noncommutative product on ${\B}$ is inherited 
from the Moyal--Vey product $\circ$ for $\W_B$ 
and the noncommutative product $\ast$ on the auxiliary bundle ${\U}$.
Let us denote the product on  ${\B}$ also by $\circ$.
The noncommutative product on  ${\U}$
contains terms in different degrees in $\hbar$, but  does not have any 
other degree bearing terms. Hence, 
the product on ${\B}$ is not preserving the grading
anymore. 

The Lie algebra structure is defined as well. For $f, g \in \Gamma(B, \B):$
\[
 [f,g] = [f,g]_\circ + [f,g]_\ast.
\]
The symplectic connection on $B$ and the connection on the auxiliary
bundle give rise to a connection on $\B$:
\[
  \nabla = \nabla^B \otimes 1 + 1 \otimes \nabla^{\F}. 
 \]
 In what
follows we will omit the tensor signs, this should not cause a
confusion. 

The curvature of this connection is 
\[
R = R^B + R^{\F}
 \]
 Following the general scheme  we want to flatten the connection
$\nabla$ by adding operators of degree $-1$ and higher.

\noindent {\bf Main theorem.}
 \begin{theorem}
\label{th:main1}
The equation $(D^{\F})^2=0$ for the connection 
\[
D^{\F}: {\cal A}^q (B,{\B}) \to  {\cal A}^{q+1} (B,{\B})
\]
has  a  solution in the form:
\begin{equation}
D^{\F}= d + \disph[{\gamma}, \cdot],
		\label{eq:flatD}
\end{equation}
 where ${\gamma} \in Inn {\B}$ is a sum of operators of degree
 $\geq -1$ , so that $d + \disph[{\gamma}, \cdot] = \nabla + \delta + r$,in
particular \ $deg (\delta) = -1$ and \[
r = \sum_{k\geq 1}r_k, \quad r_k \in gr_k \B.
\]
When $r$ satisfies the normalization condition $\delta^{-1} r = 0$
the solution is unique.

Flat sections of this connection 
are in one--to--one correspondence 
with  sections of the auxiliary bundle, $\Gamma (B, \U)$:
\begin{equation}
\label{eq:corr}
Q^\F: \Gamma (B, \U) \to \Gamma_{\text{flat}} (B, \W).
\end{equation}
\end{theorem}
%%%%%%%%%%%%%%%%%%%%%%%%%%

\begin{pf}
Main equation is $\omega + 
d \gamma + \disph[{\gamma}, \gamma] =0.$ It can be rewritten
  \[
 R + \frac{1}{2} [\nabla, \delta] + 
\frac{1}{2}[\nabla, r] + \delta^2 +\frac{1}{2} [\delta, r] + r^2 = 0. 
	\]
 For each $k > 1$ we get an equation expressing $r_k$ in
terms of $r_i, i \leq k:$
 \begin{equation}
[\delta, r_k]= [\nabla, r_{k-1}] + \sum_{i=1}^{k-2}[r_i, r_{k-i}].
\label{eq:D}
 \end{equation} 
However the equation on $r-1$ which must kill the curvature 
 gives an unusual term in $-1$  degree.
Namely, the equation on $r_1$: $\delta r_1 = R^B + R^{\F}$
gives  a solution which is a sum of two terms in degrees $1$ and $-1$ (sic!).
Indeed, in coordinates:
 \[
 r_1 = \delta^{-1}(R^B + R^{\F}) = \delta^{-1}(R^B_{ijkl}e^ie^j +
 R^{\F}_{kl})dx^k dx^l = (R^B_{ijkl}e^ie^j e^k +
 R^{\F}_{kl} e^k) dx^l.
 \]
 The operator $R^{\F}_{kl} e^k dx^l =\disph ad \{ H_{T_{kl}}e^k dx^l \}$ 
acts on a section ${\mathbf s} \in \Gamma(B, {\B})$ in the following way:
 \[
\disph  {[H_{T_{kl}} e^k dx^l, {\mathbf s}] 
=\disph  [e^k, {\mathbf s}]_\circ  H_{T_{kl}}dx^l +  \{H_{T_{kl}},
 {\mathbf s}\}e^k dx^l }
\]
where $ T_{kl}$ is an element in $\goth g$ on $M$ and hence its action 
on $\A(F)$ is defined, so it is also defined on sections of the bundle
${\B}: \quad [H_{T_{kl}}, {\mathbf s}]_\ast
 = \{H_{T_{kl}}, {\mathbf s}\} = T_{kl}{\mathbf s}.$

%%% (Also $ [H_\lambda,{\mathbf s}]_\ast = \{H_\lambda, {\mathbf s}\} $ )

The operator $\disph \text{ad}_\circ \{e^k \}H_{T_{kl}} dx^l$ is of degree
$-1$.  This term has to be added to $\delta$, the initial
$-1$--degree operator, but as soon as it is present it changes all the
equations (\ref{eq:D}) since there is not only $\delta$ in degree $-1$
anymore.

However if we change $\delta$ the iteration method can still be applied
yielding a solution for the flat connection.  In the central extension
the curvature of this connection starts from the term of 
degree
 $-2:  
\quad  {^{\F} \!\delta}^2 = {\disph}  \text{ad}_{\circ}\{(\omega_{kl} +
 H_{T_{kl}} )dx^k\wedge dx^l \}.$
 We are looking for 
the solution in the form  $^{\F} \! \delta = 
\disph A_{ik}e^i dx^k$. Applied twice it 
must give this central element:
 \[
\disph(\omega_{kl} + H_{T_{kl}}) dx^k \wedge dx^l 
=   gr_{(-2)}[\disph A_{ik}e^i dx^k, \disph  A_{jl}e^j dx^l].
 \]
 Obviously,  $A_{jl} \in \Gamma(B, \U)$  should  be 
in the form of some series in  $\omega$ and $H_T$:
 \[
A_{ik}\omega^{ij}A_{jl}= \omega_{kl} + H_{T_{kl}},  \quad \text{or} \ 
A_{ik} = \omega_{im} (1_{mk} + 1/2 \ \omega^{mn} H_{T_{nk}} + \ldots)
= \omega_{im} \left( \sqrt{ 1 + \omega^{-1}H_{T}} \ \right) ^{m} _{k}
 \]
 We have used here the skew symmetry of the forms $\omega_{kl}= -
\omega_{lk}$ and $H_{T_{kl}} = - H_{T_{lk}}$.  This series converges if
the ratio of $\omega$ and $H_{T}$ is much bigger than $1$ (i.e. size of
the fibres is very small in comparison to the base).  It means that the
fiberwise symplectic form should be much smaller than the one on the
base.

\noindent This forces us to introduce a small parameter $\epsilon$:
 \[
 (\omega_{kl} + H_{T_{kl}}) dx^k \wedge dx^l \to  ( {\omega}_{kl}
+ {\epsilon^2} H_{T_{kl}}) dx^k  \wedge dx^l.
\]
 Then the $(-1)$--degree term in the flat connection should be
\[
{^{\F} \!\delta} =  
{\disph}\omega_{km}
\left( \sqrt{ 1 + {\epsilon^2} \omega^{-1} H_{T}} \ 
\right) ^{m} _{l}dx^l \text{ad}_\circ \{e^k\}.
\]
Now the equation for the flat connection $D$ should start 
with ${^{\F} \!\delta}$
instead of $\delta$:
\[
D = {^{\F} \!\delta} + \nabla + r.
\]
Let $X$ be  the following  matrix
\[
X = \sqrt{ 1 + {\epsilon^2} \omega^{-1} H_{T}}.
\]
Let us write  how  ${^{\F} \!\delta}$ acts on 
 on a monomial
 \[
{\mathbf s} = e^{i_1}\otimes...
\otimes e^{i_m} \otimes f \otimes dx^{j_1} \wedge ...
\wedge dx^{j_n} \in gr_m (\W_B \otimes_{\C[[\hbar]]} {\U}) \otimes 
\Lambda^n {\mathcal T}^*B, \quad f \in gr_0({\U})
 \]
\begin{eqnarray*}
\lefteqn{
^{\F} \! \delta: \ e^{i_1}\otimes...\otimes e^{i_m}
\otimes dx^{j_1} \wedge ...\wedge dx^{j_n}
\mapsto}\\
& &\sum _{k=1}^m
e^{i_1}\otimes...\widehat{e^{i_k}}...\otimes e^{i_m}\otimes f \otimes  
X^{i_k}_l dx^l \wedge dx^{j_1} \wedge ...\wedge dx^{j_n}
\end{eqnarray*}
 In order to get rid of this correction term $X$ we can ``rescale''
the bundle $\W_B $ by changing 
\[
e^k \to \hat{e}^k = (X^{-1} e)^k \in \W_B \otimes_{\C[[\hbar]]} {\U}
\]
After some calculations we arrive at the following commutation relation:
 \begin{equation}
\label{eq:ura}
[\hat {e}^k,\hat {e}^l]_\circ = -i \hbar (\omega + \epsilon^2 H_T)^{kl}.
\end{equation}
On monomials in $\hat {e}, \ ^{\F} \! \delta $ has the same action as 
$\delta $ on monomials in $e$ (Lemma \ref{l:delta}).
 We can define $\ ^{\F} \! \delta^{-1} $
and the whole set up in the twisted Weyl bundle becomes 
a familiar data for Fedosov quantization of a symplectic manifold.
This way we reduce the problem to the usual Fedosov quantization and 
hence our theorem is proved.
 
We construct a flat connection on the twisted bundle by iterations, $r_i$
being monomials in $\hat {e}$.
 Uniqueness follows from the condition
$\delta^{-1}r = 0,$
which gives the same restriction as 
${^{\F} \!\delta}^{-1}r = 0.$
 The flat sections of this new connection corresponding to the sections
of the auxiliary bundle can also be obtained by the recursive procedure.
 
The only difference we must emphasize is in the characteristic class 
of deformation.
$^{\F} \!\delta^2$ gives that central element which determines the
characteristic class of deformation:
$\Omega=\omega + \epsilon^2 H_T.$
 \end{pf}
%\hspace{-2cm}{\tt char class}

%%%%%%%%%%%%%%%%%%%%%%%%

%\subsection{Weak coupling limit. MM calculus.\\}
\subsection{Quantization of the total space.\\}
\label{sec:qts}
 Here we want to prove
that the deformation we got in the previous section is indeed a
deformation of functions on $M$. Since by quantization of functions on
$M$ we understand an isomorphism to flat sections of a flat connection 
of Weyl algebra bundle over $M$ we reformulate the problem as follows.
 Namely, there is a homomorphism of noncommutative 
algebras $\alpha$ such that 
the following diagram commutes:
\[
   \dgARROWLENGTH=0.6\dgARROWLENGTH
   \begin{diagram}                 
 \node{ C^\infty(M)[[ \hbar]]} \arrow{s,l}{Q}	 
\arrow{e}{}    
  \node{\Gamma(B, {\cal F})}
\arrow{s,r}{Q^\F}
	\\
 \node{ \Gamma_{\text{\it flat}} (M, \W_M)} \arrow{e,t,..}{\alpha}    
  \node{ \Gamma_{\text{\it flat}} (B, \B)}
\end{diagram}
\]
 First line is an isomorphism of commutative algebras (\ref{eq:C(M)}).
Left vertical arrow is a quantization of a symplectic manifold $M$
(\ref{eq:Q}). 
 $Q^\F$ is a lift of the quantization $Q^{\text{fibre}}$
(\ref{eq:Qf}) to $\Gamma(B, \B)$ followed by quantization in 
the twisted Weyl algebra bundle (\ref{eq:corr}) from
Theorem \ref{th:main1}. 
  
 The Weyl algebra structure is defined on
$\Gamma (M, W_M)$ from the symplectic form  (\ref{eq:formM}) on $M$,
and 
$\Gamma (B, W_B)$ from the symplectic form $\omega_B$ on $B$. However 
\begin{proposition}
 There exists an isomorphism of {\em algebra} bundles:
\[
\beta: \Gamma(M,{\cal W}_M)  \to 
\Gamma(B,{\cal W}_B \otimes_{\C [[\hbar]]}
 \Gamma(M/B,{\cal W}_{M/B})).
\]
Let  $ D^F$ be a flat connection on $\Gamma(M/B,{\cal W}_{M/B})$, such that 
$\Gamma_{\text{flat}}(M/B,{\cal W}_{M/B}) = \U$.

Then there exists a flat connection on ${\cal W}_M$:
\[
 D^M = \bar{D} + D^F,
\]
such that its flat sections under isomorphism $\beta$ are mapped 
to the flat sections of $D^\F$ in $\B = {\cal W}_B \otimes_{\C [[\hbar]]}\U$.
It leads to the homomorphism of algebras:
\[
{\A} (M) \cong \Gamma_{\text{flat}} (B, {\B}).
\]
\end{proposition}

\begin{pf}
Consider the  local frames introduced at the end of Section 
\ref{loc}:
 $\{f_\alpha \}$
of vertical tangent bundle ${\mathcal V}M$ and  $\{e_i\}$ in 
${\mathcal T}B$ at
a point $b  = \pi( x)$ of $B$, with dual frames $\{f^{\alpha}\}$ and 
$\{e^i\}$.

\noindent  Using the connection we obtain a local frame on the
tangent bundle ${\mathcal T}M = \pi^*{\mathcal T}B \oplus {\mathcal V}M$ 
at a point $x$: 
$\{\pi^*e_i, f_\alpha \}.$ 
Let us denote the dual frame as 
$\{ \bar{e} ^i, f^\alpha \}.$ 
 
\noindent Then the form on $M$ can be written as a block matrix
(\ref{eq:formM}), so that in this frame
 \[
[f^{\alpha},f^{\beta}] 
 = {\epsilon^{2}}\sigma_b^{{\alpha}{\beta}}, \quad 
[ \bar{e} ^k, \bar{e} ^l ]  
=  ({\pi ^* \omega ^B } + {\epsilon^2} H_T )^{kl}.
\]
 The flat connection $D^F$ defines  an isomorphism of algebras 
of fibrewise deformed
functions. The map  $ \W_M \to \B: \bar{e} \mapsto \hat{e}$
defines the algebra homomorphism:
\[
\Gamma_{\text{\it fibrewise flat}}(M, \W_M)
\cong \Gamma(B, \W_B \otimes _{\C [[\hbar]]}{\U}).
\]
Indeed although the Moyal--Vey products are different, the 
product $\hat{e}$ in $\B$ and $\bar{e}$ in $ \W_M$ give the same 
$\ast$--product
(see (\ref{eq:ura})), so it gives a homomorphism of algebras.

The last step is to construct a connection 
 $\bar{D}$ in $\Gamma_{\text{\it fibrewise flat}}(M, \W_M)$ so that 
 flat sections of $D^\F$ (\ref{eq:flatD}) correspond to the flat sections of 
 $\bar{D}$. We do it just by writing for $\bar{D}$ the same formula as
for $D^\F$ but with $\bar{e}$ instead of $\hat{e}$.

The connection $\bar{D}$ together with the Fedosov connection along the
fibres $D^F$ gives rise to a single connection $D^M$ on $\mathcal W$:
\[
D^M = D^F + \bar{D}.
\]
Its flat sections $s \in \Gamma(M, {\mathcal W})$ satisfy  simultaneously
two equations 
$D^F s = 0, \quad \bar {D} s = 0.$ 
 So we get a
map from the quantization of the total space to the quantization on the
base with values in the auxiliary bundle $\U$. 
 \end{pf}

%%%%%%%%%%%%%%%%%%%%%%%%%%%%%%
\noindent {\bf Interpretation in  MM--calculus.} 
 In order to see how this homomorphism works 
we  want to introduce calculus similar to the one in \cite{MM}.

There is  some analogy with Riemannian fibrations (see for example \cite{B}). 
Mazzeo and Melrose (\cite{MM}) gave an interpretation of Hodge--Leray
spectral sequence from an analytic point of view. In particular they
introduced language similar to $b$--calculus for description of
Riemannian fibrations. The idea was to put in a small parameter
$\epsilon$, so that all horizontal differential forms had the parameter
in some degree.  So everything which came from the base was ``marked''
by this small parameter. This gave the description of terms in 
the spectral sequence by the coefficients in Taylor decomposition
with respect to
this small parameter.

Symplectic fibrations provide somewhat similar picture. One even has
the parameter naturally coming in the construction of a symplectic form
on the total space.  Indeed, when the parameter is zero one  gets a
fibrewise noncommutative product while along the base it is
commutative. The $\ast$--product then is a bidifferential expression
only in vertical coordinates. The term at the first degree of the
parameter gives a Poisson bracket in the horizontal
direction, it is a first order bidifferential expression in the 
horizontal direction.  If the fibration is trivial these bidifferential
expressions are fibrewise constant.

Given a connection (\ref{eq:split}) on $M$: ${\mathcal T}\!M 
= {\mathcal H}M \oplus {\mathcal V}\!M$ one implements the 
splitting into the structure of the product manifold 
$X = M \times [0,\epsilon_0)$, where $\epsilon_0 \ll 1$ is some fixed 
small number. (We want it to be small enough so that 
$\epsilon$ involved in 
the symplectic form $(\ref{eq:form})$ is bigger than $\epsilon_0$.)
The product $X = M \times [0,\epsilon)$ has an induced fibration,
with leaves $F$ and the base $B \times [0,\epsilon)$. 
Consider the space $\mathcal L$ of smooth vector fields on $X$ which 
are tangent to the fibres, $M$, of the product structure and  which are 
also tangent to the fibres of the fibration $M\to B$,  above $M_0 = \{\nu=0\}$.
In local coordinates ${x^j, \xi^k}$ in $M$, where the $x$'s 
give coordinates in $B$, the elements of $\mathcal L$ are the vector 
fields of the form 
\[
\sum_{j=1}^{2p}a (x, \xi, \nu)\nu \partial_{x_j}+ 
\sum_{k=1}^{2q}b(x, \xi, \nu) \partial_{\xi_k}.
\]
Consider a vector bundle $^\nu {\mathcal T}M$ for which 
$\mathcal L$ is the set of  sections
\[
{\mathcal L}= C^\infty (X,  ^\nu {\mathcal T}M).
\]
There is a natural bundle map 
$\imath_\nu: \  ^\nu{\mathcal T}M \to {\mathcal T}_XM $,
the lift of  ${\mathcal T}M $ to $X$. It is an 
isomorphism except over $M_0$, where its range is equal to ${\mathcal V}M$.
It is important to define the dual map
\[
\imath^\nu: \ {\mathcal T}^*_XM \to ^\nu\!{\mathcal T}^*M, 
\]
which range over $M_0$ is a subbundle which is naturally isomorphic to the 
bundle of forms on fibres.

Given a connection (\ref{eq:split}) on $M$: ${\mathcal T}\!M 
= {\mathcal H}M \oplus {\mathcal V}\!M$, 
the restriction of $^\nu {\mathcal T}^*M$ to the boundary, 
$M_0$, of $X$ naturally splits
\[
^\nu {\mathcal T}^*\!_{M_0}M = 
\nu^{-1} {\mathcal T}^*B \oplus   {\mathcal V}^*M,  
\quad u_0 =
\nu^{-1} \pi^* \beta +  \imath_\nu (\alpha).
\]
The exterior powers  of  $^\nu\!{\mathcal T}^*M$   
also split at the boundary so one can 
define a new bundle of rescaled differential forms on $X$. Namely,
\[
^\nu {\cal A}^k_x(M) = \sum_ {j=0}^k {\cal A}^j_x(M/B) \oplus  \nu ^{-(k-j)} 
{\cal A}^{k-j}_{\pi(x)}(B ).
\]

The symmetric powers of $^\nu\!{\mathcal T}^*M$ have the 
following decomposition
\[
^\nu S^k {\mathcal T^*}M = \sum_ {j=0}^k S^j(M/B) 
\oplus  \nu^{-(k-j)} S^{k-j}{\mathcal T^*}B.
\]
This way one can define the rescaled Weyl algebra bundle.
At a point $x \in M$ as before one can introduce a local frame $\{f_k \}$
of vertical tangent bundle ${\mathcal V}M$ 
corresponding to $\partial_{\xi_k}$ and a local frame $\{e_i\}$ in 
${\mathcal T}B$, corresponding to $\partial_{x_j}$ at
a point $b  = \pi( x)$ of $B$, with dual frames $\{f^{k}\}$ and 
$\{e^i\}$. Using the connection we obtain a local frame on the
tangent bundle ${\mathcal T}M = \pi^*{\mathcal T}B \oplus {\mathcal V}M$ 
at a point $x$.
The differential on the bundle of forms $^\nu {\cal A} ^k (M)$ is given by
\[
d = \nu \frac{dx^j}{\nu} \partial_{x_j} + d \xi^k \partial_{\xi_k}, 
\]
so that $d : ^\nu {\cal A} ^k (M) \to ^\nu {\cal A} ^{k+1} (M)$.
Similarly we can find a  
symplectic connection on  $^\nu\!{\mathcal T}M$ from the 
 initial connection on
${\mathcal T}M$. 

Then $^\nu D$, the Fedosov connection on rescaled Weyl 
algebra bundle $^\nu{\cal W}_M$
will have the Taylor decomposition into degrees of $\nu$. Since Fedosov 
connection is flat it should give an equation in each degree of $\nu$. 
The $\nu$--decomposition of $(^\nu D)^2$ must  give $0$ in 
each degree of $\nu$.
The result for quantization can be stated as follows:
\begin{enumerate}
\item
The quantization of $M$ for $\nu = 0$ is 
$C^{\infty}(B,{\A}(M/B))$
\item
The term at the first power of $\nu$ is the product on the base given by the 
Poisson bracket 
with values in the quantization of the fiber.
\item 
$n$--th power of $\nu$ allows one to write a product on the base 
with values in the product of the fibre up-to the $n$--th power in $\hbar$.
\end{enumerate}
The total space of a symplectic fibration $M \to B$
together with  
the rescaled  symplectic form and a connection preserving it   make 
such an $F$--manifold which is easily associated with an
$F$--bundle and hence gives a map from quantization of the total 
space to  a quantization of the base with values in the auxiliary bundle
of fibrewise quantizations.

%Where $\nu$ is a small number.
%\[
%x^{-1} = \{ (\frac{\omega^B}{\epsilon} + h_T) 
%(\omega^B)^{-1} \} ^{ \frac{1}{2}}
%\]
%Notice that $x$ is a matrix of the dimension of the base with entries being 
%functions on the base with coefficients being functions  on the fiber.
%Using this rescaling we can construct a symplectic connection on $M$ from 
%those on the base and on the fiber, such that the restriction 
%to the fiber will be 
%a connection on that fiber. 

%%%%%%%%%%%%%%%%%%%%%%%%%%%%%%%%%%%%%%%%%%%%%%%%%%%%%%%%%%%%

\section{Examples of symplectic fibrations and their quantization.}

 Fedosov quantization provides a way to construct a $\ast$--product on a
symplectic manifold. In the previous section we showed how Fedosov
quantization works for any symplectic fibration. However, step--by--step
calculations become complicated very quickly and explicit formulas are
readily available only in a few particular cases.  Mostly the results
which we are getting are of the type that
"under quantization  some of the variables
behave in a certain way". 

 For the trivial symplectic fibration, that is for the direct product of
two $F$--manifolds $(B, \omega^B, \nabla^B)$ and $(F, \sigma, \nabla^F)$
 we get the direct product of quantizations:
\[
\A (B \times F) = \A (B) \times _{\C[[\hbar]]} \A (F),
 \]
 with the characteristic class $ \omega^B + \epsilon \sigma $, where
$\epsilon$ can be arbitrary nonzero number (since automatically for
nondegenerate $\omega$ and $\sigma$ for any nonzero $\epsilon$ the sum is
nondegenerate in the absence of the curvature term). 

An example of a symplectic fibration with fibres 
being ${\R}^{2n}$ is considered in \cite{F2}.

 Good source of examples of symplectic fibrations 
are  cotangent bundles to fibre bundles with connections 
(see again \cite{GLS}).
 Indeed, consider a bundle $X $ over $ B$, such that $X$ is an
associated bundle to a principal $G$--bundle $P$ with a fibre $F:$ $ X =
P \times_G F.$ Then $G$ acts by symplectomorphisms on the cotangent
bundle ${\mathcal T}^*\!F$
 equipped with the canonical symplectic form $d \alpha^F$.  We
construct a new bundle $M = {\mathcal T}^*\!X$ over ${\mathcal T}^*\!B$,
 which is the pullback of $P
\times_G {\mathcal T}^*\!F$ to ${\mathcal T}^*\!B$.
 
The connection, that is a splitting of ${\mathcal T}\!M$, 
is inherited from the splitting of ${\mathcal T}X$.
Namely, ${\mathcal H}\!M$ is defined as the preimage 
of   ${\mathcal H}\!X$.

 The symplectic form on $M$ is the canonical symplectic form on the
cotangent bundle. However due to the fibre bundle structure we can 
split it into two parts -- one along the fibres and the other one 
coming from the base. There is no need to introduce the small parameter 
since we know that the form is nondegenerate as is.

 The algebraic index theorem on a cotangent bundle coincides with the
analytical index theorem of Atiyah--Singer on the initial space
\cite{NT1}.  Our example shows that this way we may obtain the relation
on the index of the fibre, the base, and the total space for any 
fibre bundle (if the  transformation group is a compact Lie group).

 Now let us move to  other  particular 
examples. Consider the following situation:
the associated bundle $M \to B$ to the principal 
$G$-bundle with a fibre being ${\mathcal T}^* G$ -- 
the cotangent bundle to $G$. This example is the inverse of the 
quantum reduction problem and it was discussed in the recent 
article \cite{F3}, which is the generalization of \cite{F1}.

 The article \cite{F1} treats a symplectic fibration 
$M \to B$ with  a symplectic fibre being 
a cylinder. 
The fibre can be
represented as  ${\C}^*= {\mathcal T}^* S^1$ -- 
the cotangent bundle to
the circle. Locally a point $z$ in $M$ can be described by 
coordinates $(x, r, \theta)$, $\theta$ being an angle coordinate on 
the cylinder and $r$ the height, while
 $x=(x_1, ..., x_{2n})$ denotes  coordinates of the base point 
${\pi(z)}$.
Then $M= P \times_{U(1)} {\C}^*$, where $P$ is a principal ${U(1)}$--bundle.
The symplectic form can be  constructed from Proposition (\ref{prop:form}). 
 Let $\lambda$ be a connection one--form on $P$. In this particular case:
\[
\lambda: \! {\mathcal T}P \to {\goth g} = \R.
\]
  Hamiltonian of
$\lambda$ is a function in $r$ only, it does not involve
$\theta$. Quantization on the fibres is the Weyl quantization like in
${\R}^{2n}$ (see (\ref{eq:R2n})). The symplectic connection is flat and
one gets from this fibrewise quantization that fibrewise flat sections
are expressions in 
 $(r+f^1), (\theta + f^2)$, where $(f^1, f^2)$ are generators 
in the Weyl algebra bundle  corresponding to $(dr, d\theta)$ 
in $T^*F$. The resulting flat connection 
in the Weyl algebra bundle on $M$ is a series in  $(r+f^1)$ 
(it does not involve the other coordinate $(\theta + f^2)$). 

Examples with  two--dimensional fibres are provided by 
any Riemannian 
surfaces (except a torus) as fibres (see \cite{McDS}).

%%%%%%%%%%%%%%%%%%%%%%%%%%%%%%%%%%%%%%%

As a nontrivial example of symplectic fibrations we consider
$S^2$--bundles.  $S^2$--bundle $\pi: M \to B$ can be considered as an
associated bundle to a principal $U(1)$--bundle, $P$. We construct the
symplectic form on $M = P \times_{U(1)}S^2$. The manifold $M$ is
presented as the symplectic reduction of $(W, \Omega)$ at a regular
value $0$ of the moment map.  The algebra of the group $U(1)$ is simply
${\R}$.
 Let ${\mathcal V}\!P \subset {\mathcal T}\!P$ be the
bundle of vertical tangent vectors. The fibre at a point $p$ being
${\mathcal V}_p\!P \subset {\mathcal T}_p\!P$.   
A connection $1$--form, $\lambda:{\mathcal T}\!_pP \to {\R}$ 
determines a horizontal subbundle by
 \[ 
{\mathcal H}_p = \{ v \in {\mathcal T}_pP \mid \lambda_p (v) =0\}
\]
This horizontal subbundle induces an injection:
$\imath_\lambda: {\mathcal V}^*\!P \hookrightarrow {\mathcal T}^*\!P,$
namely, a vertical cotangent vector $\xi \in {\mathcal V}^*_pP$ is a 
linear functional
on ${\mathcal T}_p\!P$ which vanishes on the horizontal 
subspace ${\mathcal H}_pP$. The subbundle ${\mathcal V}^*\!P$ 
inherits the standard symplectic form from ${\mathcal T}_*P \ --   \ 
d \alpha_{\mbox{can}}.$
The manifold (Weinstein universal space)
\[
W = {\mathcal V}^*\!P \times S^2
\]
carries a natural symplectic form $\Omega = pr^*_B \omega_B + 
\imath^*_{\lambda} d \alpha_{can} + pr^*_S \sigma$,
where $pr_B: W \rightarrow B$ and $pr_S: W \rightarrow S^2$ are 
the obvious projections and $\sigma$ is a 
$U(1)$--invariant volume form on $S^2$.
$\Omega$ is invariant under the diagonal action of $U(1)$. 
${\mathcal V}^*\!P$ is equivariantly diffeomorphic to $P \times \R$, 
so the moment map 
$\mu: W = P \times {\R} \times S^2 \to {\R}$ is given by 
\[
\mu(p, \eta, z) = h(z) - \eta
\]
where $h: S^2 \to \R$ is the  height function,  it is 
a moment map  for 
the action of $U(1)$ on $S^2$ by rotating about the vertical axis. 
$\eta $ is a projection $P \times \R \to \R.$ 

The level set $\mu^{-1} (0)$ can be identified with the manifold 
$P\times S^2$
by the map which takes the form $\Omega$ to 
$ pr^*_B \omega_B - d(H_\lambda)+ pr^*_S \sigma$
on $P\times S^2$, where $H(p,z) = h(z)$ is the height function on $S^2$.
This form is equivariant under the $U(1)$ action,
 to make sure it is nondegenerate  we need to introduce a small number 
$\epsilon$, so that the form on $M$ becomes 
\[
\omega =  pr^*_B \omega_B - \epsilon^2 \{ d(H_\lambda)+ pr^*_S \sigma \}.
\]
Let us  consider local coordinates at a point $m$ in $M$: $(x, \xi, \theta)$,
where $x=(x_1, ..., x_{2n})$ denotes  coordinates of the base point 
${\pi(m)}$ while  $\xi$ and $\theta$ are cylindrical 
polar coordinates in the fibre, $\xi$
gives a height function and $\theta$ is an angle. Then 
the symplectic form on the fibre is 
\[
\sigma = d\xi \wedge d\theta.
\]
The vertical vector field $\frac{\partial}{\partial \theta}$  
has a Hamiltonian
$\xi$:
\[
H_{\frac{\partial}{\partial \theta}} = \xi.
\]
The auxiliary bundle $\U$ is a bundle  of quantized
functions on fibres, it  is associated to $U(1)$--bundle $P.$ 
The connection on $\U$ is inherited from a connection $1$--form
by Proposition(\ref{th:con}). In coordinates it is as follows:
\[
\nabla^{\F} {\mathbf s} = d {\mathbf s}+ \lambda \{ \xi,  {\mathbf s} \},
\]
where $\lambda$ is a local $1$-form on the base and $\{ \cdot , \cdot \}$
is a fibrewise Poisson bracket, ${\mathbf s}$ some section of $\U$.
The curvature of this connection is:
\[
R^{\F}=\disph \text{ad} H_T = T \{ \xi,  \cdot \},
\]
where $T$ is a $2$--form on the base, the curvature of the 
connection $\lambda$. Since a Hamiltonian of any vector field is a 
linear function of just one coordinate $\xi$.
A flat connection is as follows:
\[
D^{\F} = \check{\delta} + \nabla  + r,
\]
where $\check{\delta} = \text{ad} {\omega_{kl}} {\hat{e}}^k dx^l)$
and  $ {\hat{e}}^k =
 (\sqrt{(1 + \epsilon^2 \omega^{-1} T \xi )} e)^k .$ 

 As a result one gets that the flat connection in the case of a sphere
bundle does not depend on the cylindrical angle coordinate $\theta$.
The characteristic class of the deformation with values in the auxiliary
bundle of quantization of the fibres is ${\omega^B}+ {\epsilon^2} T
\xi$.

%%%%%%%%%%%%%%%%%%%%%%%%%%%%%%%%%%%%%%%%%%%%%%%%%%%%%%%%%

\newpage

\appendix
\section{Calculations from diagrams.}

Fedosov quantization produces the system of recursively defined
equations. Fedosov proved that there are no obstructions to solutions.
 In what follows we are trying to show how to obtain these
equations step by step from diagrams, presenting all equations at
once. Since there are other situations when one has to solve systems of
recursive equations we also hope that our presentation might be useful
in some other calculations maybe of completely different origin.

\subsection{Fedosov connection.\\}
\label{chap:A1}

Let ${\cal A} ^p (W^i) =
 \Gamma (M,{\Lambda}^p {\mathcal T}^*M \otimes {\cal W}^i)$.
Let us represent the action of operators constituting the connection 
 \[
D = - \delta + \nabla + r_1 + r_2 + r_3 + r_4 + \ldots
\]  
 by arrows pointing in directions corresponding to their degree.

Namely, $r_k: {\cal A} ^0 (W^i) \to {\cal A} ^1 (W^{i+k})$ is drawn to
go from the point corresponding to the level $i$ in the first column to
the point in the second column $k$ rows down: $\delta $ goes up one row,
$\nabla$ is on the same level, $r_1$ goes down one level and so on.
Same operators act between first and second column.

\vspace{3mm}
\divide\dgARROWLENGTH by2
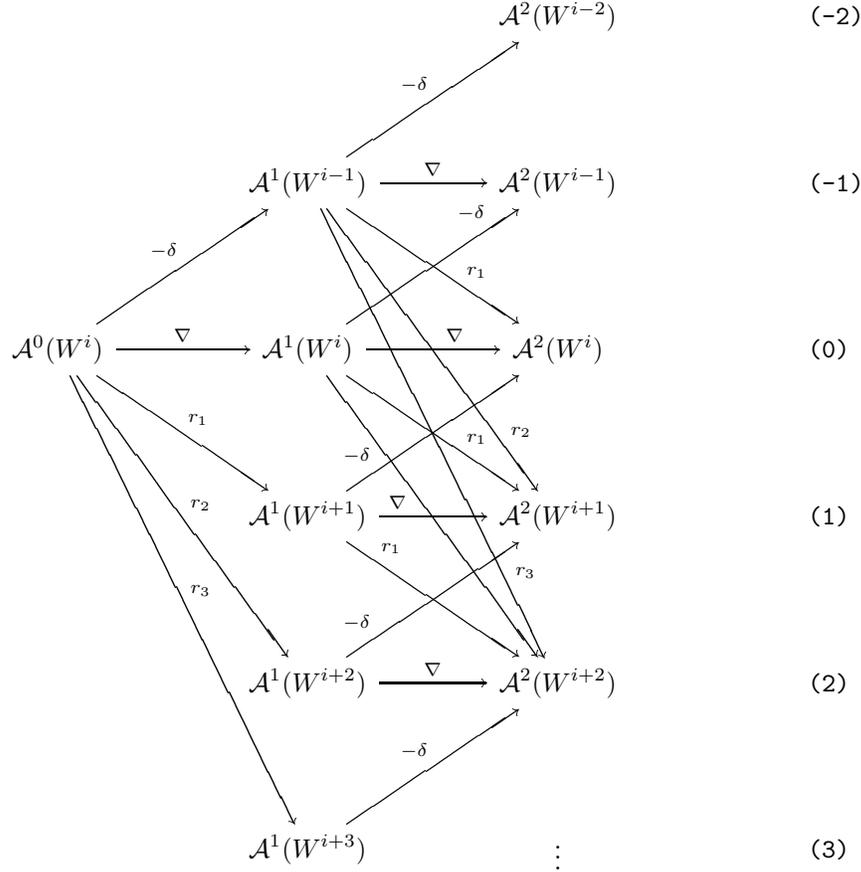
\begin{flushleft}
%
% Catcode hack to get typewriter `\' inside arg of another command
% where \verb is illegal.
\begingroup \catcode`|=0 \catcode`\\=12
   |gdef|bbb{{|tt}}%
|endgroup
%
% tighten it up a bit to fit on page in 12pt
\setlength{\dgARROWLENGTH}{1.5em}%
\noindent
$\begin{diagram}
\dgARROWLENGTH=4em
\dgARROWPARTS=6  
  \node[3]{{\cal A}^2(W^{i-2})} 
				\node{\makebox[0pt][l]{\tt\bbb (-2) }}\\
    \node[2]{{\cal A}^1 (W^{i-1})} 
	 \arrow{ne,t}{-\delta}
	\arrow{e,t}{\nabla}
	  \arrow{se,t,4}{r_1}
    	   \arrow{sse,t,5}{r_2} 
	    \arrow{ssse,t,5}{r_3}
    \node{{\cal A}^2(W^{i-1})} 
				\node{\makebox[0pt][l]{\tt\bbb (-1) }}\\
	\node{{\cal A}^0 (W^{i})} 
	    \arrow{ne,t}{-\delta}
	     \arrow{e,t}{\nabla} 
	      \arrow{se,t}{r_1}
		\arrow{sse,t}{r_2}
		 \arrow{ssse,t}{r_3}
 	\node{{\cal A}^1 (W^{i})} 
	  \arrow{ne,t,5}{-\delta}
 	   \arrow{e,t,4}{\nabla} 
	    \arrow{se,t,4}{r_1} 
	      \arrow{sse,t,5}{}
	\node{{\cal A}^2 (W^{i})}
				\node{\makebox[0pt][l]{\tt\bbb (0)  }}\\
	    \node[2]{{\cal A}^1 (W^{i+1})} 
	      \arrow{ne,t,1}{-\delta}
		\arrow{e,t,1}{\nabla}
		  \arrow{se,t,1}{r_1}
	    \node{{\cal A}^2(W^{i+1})} 
				\node{\makebox[0pt][l]{\tt\bbb (1) }}\\    
		\node[2]{{\cal A}^1 (W^{i+2})} 
		   \arrow{ne,t,1}{-\delta}
			\arrow{e,t}{\nabla}
		\node{{\cal A}^2(W^{i+2})}
				\node{\makebox[0pt][l]{\tt\bbb (2) }}\\
		      \node[2]{{\cal A}^1 (W^{i+3})} 
		          \arrow{ne,t}{-\delta} \arrow{e,!}
		      \node{\vdots}
				\node{\makebox[0pt][l]{\tt\bbb (3) }}\\ 
\end{diagram}$
\end{flushleft}

\vspace{-2cm}

 For  curvature to be equal $0$ we have to get $0$ when applying 
the connection twice to any element 
 $ a \in \Gamma (M,{\cal W})$: $D^2 a = 0.$ In other words:
the sum of arrows coming to the second column 
${\cal A} ^2 (W^i)$ should be  $0$ for every $i$.

Showing that this is true for any  element in 
 ${\cal A} ^0 (W^i) =\Gamma (M,{\cal W}^i)$ for any $i$ will do.

First two terms $\delta$ of degree $-1$  and $\nabla $ 
of degree $0$ are known,  
our purpose is to find the other terms recursively.

For every degree $i$  we get an equation on operators:

\begin{itemize}
\item Level (-2):
$\delta^2 = 0$
\item Level (-1):
$- [\delta,\nabla] = 0$
\item Level (0):
$-[\delta,r_1] + {\nabla}^2 = 0$
\item  Level (1):
$-[\delta,r_2] + [\nabla,r_1] = 0$
\item  Level (2):
$-[\delta,r_3] + [\nabla,r_2] + \frac{[r_1,r_1]}{2} = 0$
\item  Level (3):
$-[\delta,r_4] + [\nabla,r_3] + [r_1,r_2] = 0$
\end{itemize}

All these equations can be written together: $$ D^2 = ( - \delta +
\nabla + [r, .] ) ^2 = \nabla ^2 - [\delta, r] + r^2 = 0$$ where $r =
r_1 + r_2 + \ldots $ It is solved recursively: in each degree $k \geq 0
$ one gets an equation involving $r_i$ for $ i \leq k$. 

 A general proof
of the statement that there are no obstructions in this recursive
procedure is given in \cite{F} (as well as in other sources). Here let us
just show what happens in the first few equations.

{\bf Degree $-2$}. The equation is $\delta ^2 = 0$. It is satisfied 
by  the Lemma (\ref{l:delta}).

{\bf Degree $-1$}. Next one is  $[\delta,\nabla] = 0$. 
It is true by a simple calculation.

 For this equation we need that the connection $\nabla $ is torsion--free.

{\bf Degree $0$}.  Here is the first nontrivial  calculation.
We have to find such $r_1$ that
$ - [\delta,r_1] = {\nabla}^2 $.

a) Existence. First of all:
$$[\delta , \nabla ^2] =
[\delta, \nabla] \nabla - \nabla [\delta, \nabla] $$
which is $0$ by the previous equation.

There is an operator $\delta^*$ which is a homotopy for $\delta$.
\begin{equation}
	\delta^* \delta ^* = 0 , \ \ 
	\delta \delta^* + \delta^* \delta = id \ c
					\label{eq:delta}
\end{equation}
This $c$ is a number of $y$ 's and $dx$ 's, for example for
a term $y^{i_1}\ldots y^{i_p} dx ^ {j_1} \ldots dx ^ {j_q}$
this number $c = p+q $. Let us put $$r_1 =  \delta ^* \nabla ^2 $$
then indeed:
$$ [\delta, r_1] = \delta (r_1) =\delta (\delta ^* \nabla ^2) = \nabla^2$$
and also $\delta ^* r_1 = 0$.

b) Uniqueness.

Let $r'_1 = r_1 + \alpha $, such that $\delta^* \alpha = 0$. Then 
$\alpha = \delta ^* \beta$ for some $\beta$. Hence,  $\delta \delta ^* \beta 
= 0$, because $ \delta (r_1 + \alpha ) = \delta r_1$.

From (\ref{eq:delta})  we get that $\beta = \delta^* \delta \beta$
and $\alpha = \delta^* \beta = \delta ^* (\delta^* \delta \beta) = 0$

\vspace{1cm}

{\bf Degree $1$}.  $[\delta, r_2]  = [\nabla, r_1] $ gives the equation on 
the operator $r_2$.

a) Existence. Again we show
$$[\delta, [ \nabla, r_1]] = [ [\delta, \nabla], r_1] - 
[ \nabla, [ \delta, r_1]] = 0$$

b) Uniqueness. $r_2 = \delta ^* ([\nabla, r_1])$ similar to the previous one.

Recursively getting similar  equations for $r_n$
one finds the Fedosov connection. Here are first few terms:
\[
D = \delta + \nabla + \delta^{-1} \nabla^2 +
 \delta^{-1} \{\nabla, \delta^{-1} \nabla^2 \} + \ldots
\]

\subsection{Flat sections and the $\ast$--product.\\}

The Fedosov connection is flat and its flat sections are in 
one-to-one correspondence with series in $\hbar$ with 
functional coefficients. There is a map:
 \[
Q: {C}^\infty(M)[[\hbar]] \to {{\cal A}^0(W)}
 \]
Given a series:
 \[
{\mathbf a} = a_0 + \hbar a_1 + \hbar^2 a_2+ \ldots
 \]
one can find uniquely the corresponding  flat section 
of the Weyl algebra bundle
 \[
A = {\mathbf a} + A_1 + A_2 + A_3 + \ldots
 \]
 so that
 \[
\delta^{-1} (A- {\mathbf a}) = 0.
\] 
 This last condition makes the operator $Q^{-1}:{{\cal A}^0(W)} \to 
 {C}^\infty(M)[[\hbar]]$ very simple, namely it is just an evaluation of
$A \in {{\cal A}^0(W)}$ at zero value of coordinates along the fibres $W$.
This condition could be changed for any other condition fixing the 
zero section in ${{\cal A}^0(W)}$ (see \cite{EW}).

Again, the condition of flatness:
 \[ 
DA = 0
\]
 can be represented by the fact that for all $i$ sum of   operators $r_k$
which get to 
${\cal A}^1(W^{i})$ must be $0$.
It again gives a recursive system of equations.
 
We notice 
that all $r_k$ kill functions, because $r_k$ acts as adjoint operators and
functions are in the center of ${\cal A}^0 (W)$, so $r_k (a_i) = 0 $.
Hence first few  equations following from the diagram bellow are 

\begin{enumerate}
\item $ \nabla a_0 - \delta A_1 = 0$ 

\item  $ \nabla A_1 - \delta (\hbar a_1 + A_2)= 0$ 

\item $  r_1 A_1 + \nabla (\hbar a_1 +A_2) - \delta A_3   = 0$

\item $ r_2 A_1 + r_1 A_2+ \nabla A_3
				-\delta (\hbar^2 a_2 + A_4) = 0$

\item $ r_3 A_1 + r_2 A_2 + r_1 A_3
				+\nabla (\hbar^2 a_2 + A_4) - \delta A_5= 0$ 
\end{enumerate}

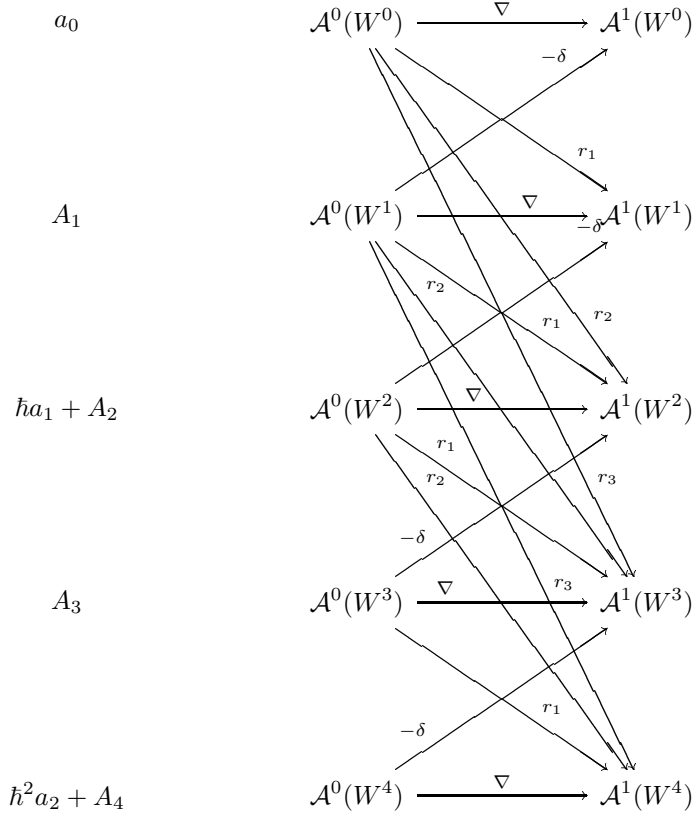
\divide\dgARROWLENGTH by2

\setlength{\dgARROWLENGTH}{1.5em}%
\noindent
$\begin{diagram}
\dgARROWLENGTH=6em
\dgARROWPARTS=6
\node{a_0} 
\arrow{e,!}
	\node{{\cal A}^0(W^{0})} 
	\arrow{e,t}{\nabla} 
	      \arrow{se,t,5}{r_1}
		\arrow{sse,t,5}{r_2}
		 \arrow{ssse,t,5}{r_3}
 	\node{{\cal A}^1 (W^{0})}\\  
\node{A_1} 
 \arrow{e,!}
	\node{{\cal A}^0 (W^{1})} 
	  \arrow{ne,t,5}{-\delta}
 	   \arrow{e,t,4}{\nabla} 
	    \arrow{se,t,4}{r_1} 
	      \arrow{sse,t,1}{r_2}
		\arrow{ssse,t,4}{r_3}
	\node{{\cal A}^1 (W^{1})}\\
\node{ \hbar a_1 +A_2} 
\arrow{e,!}
  	\node{{\cal A}^0 (W^{2})} 
	      \arrow{ne,t,6}{-\delta}
		\arrow{e,t,2}{\nabla}
		  \arrow{se,t,1}{r_1}
			\arrow{sse,t,1}{r_2}
	    \node{{\cal A}^1(W^{2})}\\ 
\node{A_3} 
	   \arrow{e,!}		
\node{{\cal A}^0 (W^{3})} 
		   \arrow{ne,t,1}{-\delta}
			\arrow{e,t,1}{\nabla}
			\arrow{se,t,4}{r_1}
		\node{{\cal A}^1(W^{3})}\\
\node{ \hbar^2 a_2 + A_4} 
	   \arrow{e,!}
		      \node{{\cal A}^0 (W^{4})} 
		          \arrow{ne,t,1}{-\delta} 
				\arrow{e,t} {\nabla}	
		      \node{{\cal A}^1 (W^4)}
\end{diagram}$

\vspace{2cm}

Let $dx^i$ be a local frame in ${\mathcal T}^*M$. Then let the 
 corresponding generators in $W$ 
be $\{e^i\}$. Then $A_i$ are of the 
form $A_{k_1 \ldots k_i}e^{k_1} \ldots e^{k_i}$.

The symplectic connection $\nabla$ locally can be written as:
\[
\nabla  = dx^j \frac{\partial}{\partial x^j} + 
\frac{i}{2 \hbar} [\Gamma_{jkl}dx^je^k e^l, \ \ ]
\]
So for the first terms of the flat section corresponding to $ { \mathbf a } = a_0$
we get:
\begin{enumerate}
\item 
$A_1 = \delta^{-1} \nabla  { \mathbf a } 
= \delta ^{-1}d{\mathbf a} = \partial_l  { \mathbf a } e^l$
\item
$A_2 = \delta^{-1} \nabla A_1 
= \delta^{-1} \nabla (\partial_l  { \mathbf a } e^l)
= \{\partial_k \partial_l  { \mathbf a } 
+ {\Gamma_{kl}}^j\partial_j  { \mathbf a } \}e^k e^l$

\end{enumerate}

The  first few terms in the $\ast$--product of two functions 
${ \mathbf a}, { \mathbf b} \in C^\infty(M)$ are:
\[
{ \mathbf a}*{ \mathbf b} = { \mathbf a }{\mathbf b }- \frac{i \hbar}{2}
 \omega^{ij} \partial_i { \mathbf a}  \partial_j {\mathbf b }
- \hbar^2  (\partial_i \partial_j {\mathbf a }+ 
\Gamma_{ij}^l \partial_l {\mathbf a }) \omega^{im} \omega^{jn}
(\partial_m \partial_n  {\mathbf b }+ 
\Gamma_{mn}^k \partial_k {\mathbf b }) + \ldots 
\]

%%%%%%%%%%%%%%%%%%%%%%%%%%%%%%%%%%%%%%%%%
\subsection{Examples of deformation quantization of symplectic manifolds.\\}

The procedure of deformation quantization requires calculations which
are not obvious and most of the time do not give nice formulas. However
in few cases one can calculate explicitly the
$\ast$--product for particular manifolds. The first trivial example 
is the quantization of ${\R}^{2n}$. Let $\{ x^1, \ldots, x^{2n} \}$ be 
a local coordinate system at some point $x \in {\R}^{2n}$.
The Darboux symplectic form in these coordinates is
\[
\omega = dx^i \wedge dx^{i+n}, \quad 1 \leq i \leq n
\]

The
standard symplectic form and the  trivial connection $\nabla =d$
 give an algebra of
pseudodifferential operators in ${\R}^{2n}$. 
Namely, using calculations from the diagrams in the last two sections 
we get the  Fedosov connection to be:
\[ 
D = d - \delta \quad \text{or in coordinates} \quad D = dx^i 
(\frac{\partial}{\partial x^i} - \frac{\partial}{\partial e^i}).
\]

Flat section of such connection corresponding to a 
function $\mathbf{a}$ under the quantization map  is as follows
\[
A = {\mathbf{a}} + e^i \frac{\partial{\mathbf a}}{\partial x^i} + 
e^i e^j \frac{{\partial^2}{\mathbf a}}{\partial x^i \partial x^j}
+ \cdots.
\]
 We see that it gives a formula for Taylor decomposition  of a function 
$\mathbf{a}$ at a point $x$. In fact the $e^i$ terms can be considered 
as jets.  Then the $\ast$--product of two flat sections 
is given by the formula (\ref{eq:Moyal}). It is easy to deduce that 
for two functions $\mathbf{a}$ and $\mathbf{b}$ the $\ast$--product is 
\begin{equation}
\begin{split}
{\mathbf a}\ast {\mathbf b} & =
\exp \  \{- i \hbar \frac{\partial}
{\partial y^k} \frac{\partial  }{\partial z^{n+k}}\}{\mathbf a (y)} 
{\mathbf b (z)} |_{y=z=x}\\
  &={\mathbf a}{\mathbf b} -
{i \hbar }\frac{\partial{\mathbf a}}{\partial x^{k}}
\frac{\partial{\mathbf b}}{\partial x^{n+k}} - \frac {\hbar^2 }{2}
(\frac{\partial^2{\mathbf a}}{\partial x^{k}\partial x^{l}})
(\frac{\partial^2{\mathbf b}}{{\partial x^{n+k}}{\partial x^{n+l}}}) + \cdots.
\end{split}
\label{eq:R2n}
\end{equation}
Let us   map $C^\infty ({\R}^{2n})$
 to differential operators 
on ${\R}^{n}$, considered as polynomials on $T^*{\R}^{n}$. Then the 
 $\ast$--product gives exactly the product of differential 
symbols.

\begin{remark}
The same scheme actually works for any cotangent  
bundle ${\mathcal T}^*M$ with the 
canonical symplectic form -- the
quantized algebra of functions on  ${\mathcal T}^*M$ is 
isomorphic to the algebra of
differential operators on $M$. This observation leads to various types of 
index theorems \cite{Gutt}, also \cite{Ezra},\cite{NT0},\cite{NT}.
\end{remark}

Explicit formulas for 
quantization of K\"ahler manifolds  were given in  \cite{K}.
Tamarkin (\cite{T}) showed that for symmetric K\"ahler 
manifolds the Fedosov
connection can be constructed to have  only three summands, it has 
 no terms of degree higher than 1, i.e.
\[
D= \delta + \nabla + r_1.
\]
This gives compact formulas for the $\ast$-product on such manifolds as 
 the 2--sphere, any projective space ${\C}P^n$,  and
Grassmanians.
The case of ${\C} P^n$ was also  considered in (\cite{Germans}).
%%%%%%%%%%%%%%%%%%%%%%%%%

\end{document}